\documentclass[review]{elsarticle}

\usepackage{lineno,hyperref}
\modulolinenumbers[5]

\journal{Journal of \LaTeX\ Templates}

\newcommand{\R}{\mathbb{R}}
\renewcommand{\exp}{\mathrm{exp}}
\newcommand{\expc}{\mathrm{expc}}

\usepackage{graphicx,amsmath,amssymb,natbib}
\usepackage{caption}
 \usepackage{easybmat}

\begin{document}

\begin{frontmatter}

\title{Computational of periodic oscillations and related bifurcations in the Hodgkin-Huxley model}

% \address{Radarweg 29, Amsterdam}
% \fntext[myfootnote]{Since 1880.}

\author[mymainaddress]{A. Balti}
\ead{aymen.balti@univ-lehavre.fr}


\author[mymainaddress]{V. Lanza\corref{mycorrespondingauthor}}
\cortext[mycorrespondingauthor]{Corresponding author}
\ead{valentina.lanza@univ-lehavre.fr}

\author[mymainaddress]{M. A. Aziz-Alaoui}

\address[mymainaddress]{Normandie Univ, France; ULH, LMAH, F-76600 Le Havre; FR CNRS 3335, ISCN, 25 rue Philippe Lebon 76600 Le Havre, France}
\ead{aziz.alaoui@univ-lehavre.fr}

% \author[mymainaddress,mysecondaryaddress]{Elsevier Inc}
% \author[mysecondaryaddress]{Global Customer Service\corref{mycorrespondingauthor}}
% \cortext[mycorrespondingauthor]{Corresponding author}
% \address[mymainaddress]{1600 John F Kennedy Boulevard, Philadelphia}
% \address[mysecondaryaddress]{360 Park Avenue South, New York}

\begin{abstract}
The Hodgkin-Huxley equations constitute one of the more realistic neuronal models in literature and the most accepted one. It is well known that,
depending on the value of the external stimuli current, it exhibits periodic solutions, both stable and unstable.\\
Our aim is to detect and characterize such periodic solutions, exploiting a robust and manageable technique, mainly based on harmonic balance method.  
\end{abstract}


\end{frontmatter}

\linenumbers

\section{Introduction}

\subsection{Biological motivation }

\label{intro} Scientists have been always fascinated by the
functioning of the human brain and always attempted to understand
its complexity.

Only in 1899 Santiago Ramon y Cajal, exploiting the experimental
techniques developed by Camillo Golgi, was able to discover that the
nervous system is made by individual cells, later called neurons \citep{Cajal}.
His studies laid the foundations for the so-called ``Neuron Doctrine''
and gave him the Nobel Prize in Physiology and Medicine in 1906
\citep{Glickstein}. Although some scientists suggest to rethink the
Neuron Doctrine \citep{Bullock}, it remains the pillar of modern
neuroscience.

It is worth observing that in each individual there is not an only
type of neuron, but several ones. Nevertheless, they share many
common properties. From the cell body (called soma) starts a number
of ramifying branches called dendrites. These structures constitute
the input pole of a neuron. From the soma originates also a long
fiber called the axon. It is considered as the output line of a
neuron since through the axons terminals, called synapses, the
exchange of information with other neurons takes place
\citep{Scott,tuckwell}.

In fact, the basic elements of the communication among the neurons
are pulsed electric signals called action potentials or spikes. The
neuronal cell is surrounded by a membrane, across with there is a
difference in electrical charge. The charge in membrane potential
depends on the flow of ions, especially Sodium (Na$^+$), Potassium
(K$^+$) and Calcium (Ca$^{++}$). If the membrane potential exceeds a
certain threshold, then the neuron generates a brief electrical
pulse, that propagates along the axon. Finally the synapses transfer
this electrical signal to the other surrounding neurons
\citep{Scott,Keener, Izhikevich}.

Moreover, a neuron can exhibit rich dynamical behaviors, such as
resting, excitable, periodic spiking, and bursting activities. In
particular, the ability of periodic firing has been recorded in
isolated neurons since the 1930s \citep{Wang} and Hodgkin
\citep{Hodgkin1948,Izhikevich} was the first to propose a
classification of neural excitability, depending on the frequency of
the action potentials generated by applying external currents.

In literature several nonlinear dynamical systems have been proposed
to suitably model the dynamics of the electrical activity observed
in single neurons.
%In Neuroscience, several mathematical models for describing the
%neurons activity have been proposed.
However, the paper by Hodgkin and Huxley \citep{HH} on the
physiology of the giant axon of the squid remains a milestone in the
science of nervous system and at present many properties of the
model proposed therein still have to be disclosed.

It is known \citep{Chaos-HH,Hassard,Rinzel} that, depending on the
value of the external current stimuli $I$, the Hodgkin-Huxley (HH)
model exhibits periodic behaviors. Moreover, for a certain range of
$I$, it shows hard oscillations \citep{lanza2,Minorsky}, that is a coexistence of a stable
equilibrium and a stable limit cycle. Thus, this implies the
existence of unstable limit cycles in order to separate the two
basins of attraction. Few authors have been able to characterize
these periodic solutions, how they emerge and disappear depending on
the intensity of the external current
\citep{Chaos-HH,Hassard,Rinzel}.

In general, it is not an easy task to detect a periodic solution of
a nonlinear dynamical system. Several methods are exploited to
predict the existence of limit cycles and to study their stability,
both in time and frequency domain \citep{Ascher,Kundert,Mickens,mees,lanza1}.
In particular, for the HH model, this is even more difficult due to
the high nonlinear structure of the system.

Since we aim to characterize and numerically approximate all the periodic solutions exhibited by the HH model,
depending on the value of the external current stimuli, let us recall briefly the structure of HH model and its main characteristics.

\subsection{The Hodgkin-Huxley model}

The Hodgkin-Huxley model for a neuron consists in a set of four
nonlinear ordinary differential equations in the four variables
$X=(V,m,h,n)$, where $V$ is the membrane potential, $m$ and $h$ are
the activation and inactivation variables of the sodium channel and
$n$ is the activation variable of the potassium current. The
corresponding equations are the following \citep{HH}:

\begin{equation}\label{eq:HH}
\left\{
\begin{array}{l}
C\dfrac{dV}{dt}  =I-[\overline{g}_{Na} m^3 h(V-E_{Na})+\overline{g}_{K} n^4(V-E_{K})+\overline{g}_{L}(V-E_{L})],\\
\\
\dfrac{dn}{dt}  = \alpha_{n}(V)(1-n)-\beta_{n}(V)n,\\
\\
\dfrac{dh}{dt}  = \alpha_{h}(V)(1-h)-\beta_{h}(V)h,\\
\\
\dfrac{dm}{dt}  = \alpha_{m}(V)(1-m)-\beta_{m}(V)m,
\end{array}
\right.
\end{equation}
%
%\begin{equation}\label{eq:HH}
%\left\{
%\begin{array}{l}
%C\frac{dV}{dt}  =\frac{1}{C_M}[I - (120m^3 h(V - 50 ) + 36n^4 (V + 77) + 0.3(V - 54.4))]\\
%\frac{dn}{dt}  = \frac{1}{10}(1-n) \expc(\frac{-55-V}{10})-\frac{1}{8}n~\exp(\frac{-V-65}{80})\\
%\frac{dh}{dt}  = (1-h)\exp(\frac{-V-65}{20})-\frac{h}{1+\exp(\frac{-V-35}{10})}\\
%\frac{dm}{dt}  = (1-m)\expc(\frac{-40-V}{10})-4m \,\,\,\exp(\frac{-V-65}{18}),
%\end{array}
%\right.
%\end{equation}
where $I$ is the external current stimulus, $C$ is membrane conductance, $\overline{g}_{i}$ are the shifted Nernst equilibrium potentials, $E_i$ are the maximal conductances, $\alpha(V)$ and $\beta(V)$ are functions of $V$, as follows:

\begin{equation*}
\begin{array}{lcl}
 \alpha_{n}(V)=0.1 \expc(0.1 (10+V) ),& & \beta_{n}(V) =\exp(V/80)/8,\\
\alpha_{h}(V)= 0.07\exp(V/20), && \beta_{h}(V) = 1/( 1+\exp( 0.1(30+V)) ),\\
\alpha_{m}(V)=\expc(0.1(25+V)),& & \beta_{m}(V) =4 \exp(V/18),
\end{array}
\end{equation*}
% $$\alpha_{n}(V)=0.1 \expc(0.1 (10+V) ),\qquad \beta_{n}(V) =\exp(V/80)/8,$$
%
% $$\alpha_{h}(V)= 0.07\exp(V/20),\qquad \beta_{h}(V) = 1/( 1+\exp( 0.1(30+V)) ),$$
%
% $$\alpha_{m}(V)=\expc(0.1(25+V)),\qquad \beta_{m}(V) =4 \exp(V/18),$$
and $\expc(x)$ is given by \citep{Hassard}  %(attention: hassard et guckenheimer n'ont pas la même fonction!!)
$$\expc(x)=\left\{
\begin{array}{l}
  \dfrac{x}{\exp(x)-1}~~~~~x\neq 0\\\\
   1~~~~~~~~~~~x=0.
\end{array}
\right.$$
Finally, the typical values for the other parameters are \citep{Hassard, Chaos-HH}:
$$E_{K}=12 \;mV,~~E_{Na}=-115\; mV,~~E_{L}=-10.599 \;mV$$
$$\overline{g}_{K}=36 \;mS/cm^2, \overline{g}_{Na}=120 \;mS/cm^2,\overline{g}_{L}=0.3 \;mS/cm^2.$$
%
%It is worth noting that the numerical values of the different
%parameters in (1) are not the common ones found in
%literature \citep{Hassard,Chaos-HH}. This set of parameters has been
%provided in the original HH paper \citep{HH}, but then the Authors
%decided to shift the membrane potential in order to have the resting
%potential at $V\approx 0$ mV. According to \citep{Izhikevich}, we
%have chosen the parameters such that the membrane potential is back
%to its natural value, that is $V\approx -65$ mV. It is possible to
%see that this choice does not change the dynamical behaviors of the
%model.

For small values of the current stimulus $I$ the system exhibits a
stable equilibrium point. If $I$ is increased, then a stable  periodic solution with large
amplitude appears, while the equilibrium point
remains stable. This means that necessarily unstable solutions are
present, in order to separate the two basins of attraction. In
\citep{Hassard}, the Author shows that, depending on the value of
$I$ the HH model presents from one to three unstable limit cycles.
Moreover, in a certain range of values of $I$ the equilibrium point
becomes unstable, but finally the stable periodic solution
disappears through an Hopf bifurcation and the equilibrium point
regains its stability.

At present, few works about the detection of the periodic solutions
of the HH model exist in literature, for example
\citep{Hassard,Chaos-HH,Rinzel}, since because of the high dimension of
the system and its high nonlinearity it is not an easy task. Furthermore, the existing works
approach the problem by exploiting different methods (finite differences,
collocation or shooting methods), that are not so simple to handle.

Our aim is to characterize and numerically approximate all the periodic solutions
exhibited by the HH model \eqref{eq:HH}, depending on the intensity of the external current stimuli $I$, through a joint application of shooting,
collocation and harmonic balance methods \citep{Ascher,mees}.
In particular, we show how the harmonic balance method permits to deduce the most complex and interesting part of the bifurcation diagram
in a more performant way with respect to the other methods.

%  We exploit a collocation method that is numerically very robust and extremely handy.
%
% Our aim is to characterize the periodic solutions of the HH model,
% exploiting a collocation (harmonic balance?) method that have been
% revealed to be numerically very robust and really handy.

The paper is structured as follows:  in
Section 2 the basics of collocation and harmonic balance methods are briefly
recalled. In Section 3 we show how it is possible to obtain the bifurcation diagram of the HH model, by exploiting shooting, collocation and harmonic balance methods.
Furthermore, all the bifurcations are analysed via the harmonic balance method and Floquet analysis. Finally, concluding remarks are offered in Section 4.

\section{Computation of periodic solutions}
%\label{sec:3}
%Le calcul pratique des solutions periodiques se base sur la resolution d'equations non lineaire en dimension finie de la forme

Classical methods, based on integration schemes, such as the shooting method, are generally sensitive to the stability properties of solutions. In this section we review two methods that overcome this problem and belong to a class of spectral methods. Indeed, we will use collocation and harmonic balance methods, whose main idea is the approximation of the exact solution by the projection on a finite-dimensional sub-space.

\subsection{Collocation methods}
\label{sec:4}

Let us consider an autonomous dynamical system
\begin{equation}\label{eq:sist}
  \dot x=f(x)
\end{equation}
where $f$ is a vector field defined on $\R^n$, $n\geq 1$, and $x\in \R^n$.
A solution $x=X$ of a continuous-time system is periodic with least
period $T$ if $X(t+T)=X(t)$ and $X(t+\tau)\neq X(t)$ for $0<\tau<T$.
This periodic solution $X$ of least finite period $T>0$ of the
system corresponds to a closed orbit $\Gamma$ in $\R^n$.
On this orbit each initial time corresponds to a location $x=x_0$.

It is interesting to notice that searching a periodic solution of an ODE is equivalent to the resolution of a boundary value problem (BVP). In fact, if $x=X$ is a $T$-periodic solution of  \eqref{eq:sist}, then it is solution of the following BVP:

\begin{equation}\label{eq:problim1}
\begin{cases}
 \dfrac{dx}{dt}=f(x)\\
   x(0)=X(0)=X(T).
\end{cases}
\end{equation}

System \eqref{eq:problim1} belongs to the general class of nonlinear
boundary value problems and several methods have been proposed in
literature to solve it \citep{Ascher}. The more intuitive
one is probably the shooting method \citep{Kuznetsov}. The idea is to find an initial
condition $X_0$ and a period $T$ such that $X(0)=X_0=X(T),$ where
the unknown is the couple $(X_0 , T)$, with $X_0 \in  \R^n$ and $T
\in \R^{+}$. If we define the functional $G$ as
$$G(X_0, T ) = \varphi(X_0, T ) - X_0,$$
where $\varphi(X_0, T )$ is the solution of the Cauchy problem $\dot x = f (x)$ with the initial condition $x(0) = X_0$, then it is straightforward to see that a periodic solution of \eqref{eq:problim1} is a zero of $G$.
Unfortunately, this method fails when the limit cycle under study is unstable, due to the errors of the numerical integration of the equations.

Since our purpose is to detect all the periodic solutions of the HH model, both stable and unstable ones, we have decided to exploit a collocation method, which is independent on the stability of the periodic solutions under consideration. We briefly recall the main properties of this method.

%\paragraph{Periodic boundary value problem }
%First of all, since $T$ is usually unknown, with a simple translation of the time scale, it is possible to transform
First of all, since $T$ is usually unknown, with a simple normalization
%$t$ by $T$, so let $\tau=\frac{t}{T}$,
of the time scale, it is possible to write \eqref{eq:problim1} on the interval $[0,1]$:

\begin{equation}\label{eq:problim2}
\begin{cases}
   \dfrac{d u}{d \tau}=Tf(u)\\
   u(0)=u(1),
\end{cases}
\end{equation}
where $\tau=\frac{t}{T}$ is the new time variable.
Clearly, a solution $u(\tau)$ of \eqref{eq:problim2} corresponds to a $T$-periodic solution of \eqref{eq:problim1}. However, the boundary condition in \eqref{eq:problim2} does not define a unique periodic solution. Indeed, any time shift of a solution to the periodic BVP \eqref{eq:problim2} is still a solution. Thus, an
additional condition has to be appended to problem \eqref{eq:problim2} in order to select a solution among all those corresponding to the cycle.

The idea of the collocation methods for BVPs is to approximate the analytical solution by a piecewise polynomial vectorial function $P(t)$ belonging to  $\R^n$ that satisfies the boundary conditions and the original problem on selected points, called collocation points.\\
Let us consider the partition $0 =t_{0} < t_{1} < \dots < t_{N}=1$, then the approximated solution $P(t)$ has the following form:
$$P(t)/_{[{t_i},{t_{i+1}}]}=P_{i}(t),\qquad i=0,\dots,N-1,$$
where $P_{i}$ is a polynomial of degree $m$ for all $i=0,\dots,N-1$, and $P_{i}(t_{i+1})=P_{i+1}(t_{i+1})$, in order to have a continuous polynomial  on the whole interval $[0,1]$. \\
On each sub-interval $[{t_i},{t_{i+1}}]$ we introduce the collocation points
$$t_{ij} = t_i + \rho_{j}(t_{i+1}-t_{i}),~i=0,~\dots,~N-1,~ j = 1,~\dots,~m,$$
where $0\leq \rho_{1}< \dots < \rho_{m} \leq 1$.

Then, the request that $P(t)$ satisfies the BVP \eqref{eq:problim2} on the collocations points leads to the following non-linear algebraic system of equations:
\begin{equation} \label{eq:colloc}
\frac{1}{T} \dot P(t_{ij})=f(P (t_{ij} )), \,\, i = 0, \dots, N-1, j = 1, \dots, m,
\end{equation}
%Here, the collocation points are given by
% $$t_{ij}=t_i + \rho_{j}(t_i-t_{i+1}),\,\,\, 0\leq \rho_{j}< \dots < \rho_{m} \leq 1.$$
%
with the boundary conditions
\begin{equation} \label{eq:cond}
P (0) = P (1).
\end{equation}
%%where $T$ is the unknown period, then
Thus, the state vector is given by
$$U=(P(0),P(t_{ij})_{0\leq i \leq N-1,1 \leq j \leq m},P(1),T)\in \R^{q},$$
where $q=m \cdot N \cdot n +2n+1,$ $m\cdot N$ is the number of unknowns in \eqref{eq:colloc}, and $n$ is the dimension of $X$ in \eqref{eq:problim2}. Moreover, a further condition is necessary to determine the unknown parameter $T$.
%Note dans notre contexte il y a beaucoup de maniere d'ajouter la derniere contrainte ... on fixe une composante... une derivee
%...contrainte par produit scalaire...la derniere est la plus regoureuse mais pas simple a implementer.

Several solvers using collocation methods have been proposed in the
literature. For example, COLSYS/COLNEW \citep{COLNEW,COLSYS} and AUTO \citep{AUTO} which use collocation methods by
using Gaussian points, or bvp4c, bvp5c and bvp6c
\citep{Shampine,JACEK} which use routines with \textit{Lobatto} points.

In our study, we will use the bvp4c one, with 3 \textit{Lobatto}
points on each subinterval. In this case $P(t)$ is a piecewise
vectorial cubic polynomial, $P(t)\in C^{1}([0,~1]^n)$
\citep{Shampine,JACEK}.

The bvp4c methods controls the error $||u(t)-P(t)||$ indirectly, by
minimizing the norm of the residue $r(t)=||\dot{P}(t)-F(P(t))||$.
Thus, $P(t)$ is considered as the exact solution of the following
perturbed problem
$$\dot{P}(t)=F(P(t))+r(t),~~H(P(0),P(1))=\delta,$$
where $H$ is a function of $P(0)$ and $P(1)$ that describes the boundary conditions, and $\delta$ is the associated residue
(in our case, $\delta=(P(1)-P(0))$).\\
The basic idea of this technique is to minimize the residue over
each sub-interval $[x_i,~x_{i+1}]$, to adjust the mesh as one goes
along, such that the norm of the residue $r$ and $\delta$ tend to zero. This
assures that the approximated solution $P(t)$ converges to the exact
solution of the problem. It is worth noting that this mesh
adaptation method is able to obtain the convergence even in case of
bad initial conditions and with an approximation of order 4, i.e.
$||u(t)-P(t)||< C h^4,$ where $C$ is a constant and $h$ is the
maximum mesh step. For further details, see \citep{Shampine}.

\subsection{Harmonic Balance (HB) method }
\label{sec:5}
It is worth noting that any periodic smooth function can be represented as an infinite Fourier series
%
%The method of harmonic balance is based on the approximation of the solution of (\ref{eq:problim1}),(it a result of Weierstrass theorem),   by the following finite Fourier series
\begin{equation}\label{eq:problim5}
 X(t)= A_0  + \sum_{k=1}^{\infty} \left(A_k \cos \left(k \dfrac{2\pi}{T} t\right) + B_k \sin \left(k \dfrac{2\pi}{T} t\right)\right),
\end{equation}
where
\begin{align}\label{eq:Four_coeff}
A_0 & =\dfrac{1}{T}\displaystyle\int_{0}^{T} X(t) dt,\nonumber\\
A_k &=\dfrac{2}{T}\displaystyle\int_{0}^{T} X(t)  \cos \left(k \dfrac{2\pi}{T} t\right)dt,~~~~~~~k=1,...,\infty.\\
B_k & =\dfrac{2}{T}\displaystyle\int_{0}^{T} X(t)  \sin \left(k \dfrac{2\pi}{T} t\right)dt,\nonumber
\end{align}

The idea of the harmonic balance method \citep{mees} is to search for an approximation of the solution of \eqref{eq:problim1}
as a truncated series

\begin{equation}\label{eq:problim6}
 \tilde{X}_{K}(t)= A_0 +\sum_{k=1}^{K} \left(A_k \cos \left(k \dfrac{2\pi}{T} t\right) + B_k \sin \left(k \dfrac{2\pi}{T} t\right)\right),
\end{equation}
where $K$ is the number of harmonics taken into account.

If the periodic solution $X(t)$ is smooth, then the truncated series $\tilde{X}_{K}(t)$ converges to $X(t)$
rapidly, without exhibiting the Gibbs phenomenon \citep{Urabe}.%and $H$ is the number of harmonic.\\

We note $\overline{X}=(A_0, A_1, B_1, ..., A_K, B_K)$ the vector of Fourier coefficients, and
$$e_{j}(t)=
  \begin{cases}
 1~~~~~~~~~~~~~~~\text{if}~~j=0,\\
\cos \left(j \dfrac{2\pi}{T} t\right)~~~\text{if}~~j=2k,~~~~~~k=1,...K.\\
\sin \left(j \dfrac{2\pi}{T} t\right)~~~\text{if}~~j=2k+1,
\end{cases}
$$
the Fourier base.
%
%We recall that $\forall t_0\in [0,~T]$, if $X(t_{0})$ continuous at $t_0$, then the approximated Fourier series $\tilde{X}(t)$ converges to the solution  $X(t_{0})$. Moreover, the convergence speed can be controlled according the following result :
%
%In non-stiff problems, the periodic solution $X(t)$ is smooth, its Fourier series does not exhibit the Gibbs phenomenon. The Fourier series $\tilde{X}_{H}(t)$ converge to $X(t)$ rapidly and uniformly. And if $X(t)$ is $C^{\infty}$, the converge is exponential order, For more detail see ...

Furthermore, if $X(t)$ is a $T$-periodic continue function, then the $n^{th}$ trigonometric Lagrange interpolation polynomial of $X(t)$ with equally spaced nodes is the following \citep{zygmund}
$$L_n(X(t),t)=\sum_{j=0}^{n} X(t_{j})l_{j}(t),$$
where
$$l_{j}(t)=\dfrac{\sin\left(\left(n+\dfrac{1}{2}\right)\left(t\dfrac{T}{2\pi}-t_j\right)\right)}{\sin\left(\dfrac{1}{2}\left(t\dfrac{T}{2\pi}-t_j\right)\right)},\qquad t_j=\dfrac{jT}{2n+1},\qquad j=0,1,...,2n.$$
Then, the functions $(l_{j})$ constitute an Hermitian base of the periodic functions space.

Let us suppose $f$ in \eqref{eq:problim1} to be a polynomial of degree $d$, and let us consider $n=d\times K$. Let $Y_{F}(\overline{X})$ be the coefficients of $f(\tilde{X}_{K}(t))$ in the Fourier base $(e_{j})_{j=0,\dots,2K}$ and $Y_{L}(\overline{X})$ the components of $f(\tilde X_{K}(t))$ in the Lagrange base $(l_{j})_{j=0,\dots,2n}$. It is easy to see that $$(Y_{L}(\overline{X}))_{j}=f(X_{K}(t_{j})).$$
Moreover, it is possible to deduce $Y_{F}(\overline{X})$:
$$ Y_{F}(\overline{X})= P\Gamma^{-1} Y_{L}(\overline{X}),$$
where $P$ is the projection matrix 
$$P=\begin{pmatrix}
     I_{2K+1} & 0 & \dots & 0
    \end{pmatrix}\in \R^{2K+1, 2n+1},$$
$I_{2K+1}$ is the identity matrix of rang $2K+1$, and $\Gamma^{-1}$ is the transition matrix between the two bases $(l_{j})_{j=0,\dots,2n}$ and $(e_{j})_{j=0,\dots,2n}$. It is easy to notice that the elements of $\Gamma^{-1}$ can be determined as the values of the functions $(e_{j}(t))_{j=0,\dots,2n}$ in the nodes $t_j$:
$$\Gamma^{-1}=\left(\begin{matrix}
 1                         &\cos(1\times t_0)  &\sin(1\times t_0)   &\cdots\cdots\cdots                   &\cos(n\times t_0)    &\sin(n\times t_0)\\
\vdots        &\vdots       &\vdots        &\cdots\cdots\cdots                   &\vdots         &\vdots\\
\vdots         &\vdots       &\vdots        &\cdots\cdots\cdots                   &\vdots         &\vdots\\
1                          &\cos(1 \times t_{2n})    &\sin(1 \times t_{2n})     &\cdots\cdots\cdots                   &\cos(n\times t_{2n})    &\sin(n \times t_{2n})\\

\end{matrix}
\right)$$

%
%Let $Y(\overline{X})$ be the Fourier coefficients  of $f(\tilde{X}(t))$. It is worth to see that in the Lagrange base the components of  $f(\tilde{X}(t))$ are easy to evaluate, and then, through  the transition matrix $\Gamma$, it is possible to return to the Fourier base and deduce Y(X) the Lagrange base, easy to evaluate, and then the transition matrix $\Gamma$,  it
%Finally, from (8) it yields
% Let $Y_{F}(\overline{X})$ be the Fourier coefficients of $f(\tilde{X}_{K}(t))$ and $Y_{L}(\overline{X})$ the components of $f(\tilde X_{K}(t))$ in the Lagrange base.
% 
% 
% $P$ is the projection operator between both sub bases fourier's of orders H and K respectively.\\
It is worth observing that the choice of $n=dK$ and the utilisation of the projection matrix $P$ permit to avoid a sort of aliasing phenomenon \citep{hesthaven}, that could take place if $n=K$. Moreover, this technique can be generalized to the case of a non-polynomial nonlinearity and, in this case, $n$ can be determined by looking at the convergence rate of the Fourier coefficients with respect to the number of considered harmonics.\\
As far as we know, this joint exploitation of Fourier and Lagrange basis in the harmonic balance method is an original approach. It is worth noting that the change of basis between the Fourier and the Lagrange ones permits to avoid to find directly the Fourier coefficients of $f(\tilde{X}_{K}(t))$ by using \eqref{eq:Four_coeff}. In particular, this is extremely useful in the case of high nonlinear systems, such as the HH model, since formulas \eqref{eq:Four_coeff} would require huge calculations, while our technique permits to obtain the Fourier coefficients in a more efficient way.

Finally, from \eqref{eq:problim6} we obtain
\begin{equation}
\tilde D\overline{X}=Y_F(\overline{X}),
\end{equation}
where\\
$\tilde D=D\otimes I_{n}$, and $D$ the matrix of differential time operator in the base $(e_j)$:

$$D=\left(\begin{matrix}
 0       &\cdots       &\cdots      &\cdots            &\cdots   &0  \\
 \vdots  &D_1          &\ddots      &\ddots            &\ddots   &\vdots  \\
 \vdots  &\ddots       &\ddots      &\ddots            &\ddots   &\vdots  \\
 \vdots  &\ddots       &\ddots      &D_k               &\ddots   &\vdots  \\
 \vdots  &\ddots       &\ddots      &\ddots            &\ddots   &\vdots  \\
 0       &\cdots       &\cdots      &\cdots            &\cdots   &D_K \\

\end{matrix}
\right)$$\\
and $D_k$ is the $2\times2$ matrix of differential time operator in the sub-base $(e_{2k} ,~ e_{2k+1}):$
 $$D_k=
 \left( \begin{matrix}
 0       &k\\
 -k     &0
\end{matrix}
\right).$$

%
%In here, our goal is determine of $Y(X)$, ie the component of $F(X)$ according to the component of X. With help Lagrange polynomial, the component of $F(X)$ is the vector $F(X(t_i))$ that is depends on $X(t_i)$,
%the transform matrix between two bases $B_1$ and $B_2$ following by
%$$\Gamma^{-1}=\left(\begin{matrix}
% 1                         &\cos(1\times t_1)  &\sin(1\times t_1)   &\cdots\cdots\cdots                   &\cos(H\times t_1)    &\sin(H\times t_1)\\
%\vdots        &\vdots       &\vdots        &\cdots\cdots\cdots                   &\vdots         &\vdots\\
%\vdots         &\vdots       &\vdots        &\cdots\cdots\cdots                   &\vdots         &\vdots\\
%1                          &\cos(1 \times t_H)    &\sin(1 \times t_H)     &\cdots\cdots\cdots                   &\cos(H\times t_H)    &\sin(H \times t_H)\\
%
%\end{matrix}
%\right)$$\\
% and the inverse transformation with help of $\Gamma$.

\section{Numerical results}
\subsection{Continuation technique and Newton's algorithm}
%\label{sec:4}

In this work, we are interested in detecting periodic solutions of HH system, depending on the values of the external current $I$. We will use a continuation method \citep{AUTO}, where at each iteration we fix the parameter $I$ and determine the unknown periodic solution $X(I)$ and its period $T(I)$, exploiting either the  collocation method or the harmonic balance method and starting from the solution found at the previous step.\\
Therefore, the numerical computation of a periodic solution is based on the resolution of the nonlinear algebraic system
\begin{equation} \label{eq:nonlinear}
F(X,T,I)=0,
\end{equation}
issue of one of the two methods presented above (collocation or HB method), where $X$ is the state vector, $T$ is the period approximation, and $I$ is the parameter.

The continuation method implementation exhibits two main
difficulties: on the one hand, the construction of a right initial
solution, and on the other hand the progression of the algorithm for critical values
of parameter $I$, that are turning points. It is possible to ride
out the first one by starting from a stable limit cycle branch, that can be
suitably numerically approximated with a fourth-order Runge-Kutta
method. For the second one, Chan and Keller \citep{KELLER} proposed the arc-length
continuation method.In our case, since in a neighborhood of the turning points the
function $T(I)$ is monotone \citep{Rinzel},  we can use $T$ as
parameter in order to follow the branch continuation.

%
%Our algorithm based in predictor-corrector (PC) method. The last method consists of the prediction of the next point. Suppose that a solution point $x_0$, and h the step associate of the unknown parameters, in our case $T$ or $I$, for sufficient choice, the branch is constant locally. and the next step is the correction for goal obtained a point in branch, in here, we use the Newton's method, end we stop if the $residue~<~\epsilon$, where $\epsilon$ present of size control. Further details on the topic can be found in [].

Let us remark that, in our case, the branches are locally linear, so for good choices of $T$ and $I$ in our algorithm the next point is obtained via a correction and by using Newton method \citep{KELLER}.

\subsection{Bifurcation diagram and branches of periodic solutions}
%
%\label{sec:5}
Both collocation and harmonic balance methods are appropriate for detecting all the periodic solutions, both the stable and unstable ones. However, the collocation method requires to adapt the mesh for giving the best results, so a huge number of nodes is needed in order to get a suitable approximation of the solution. On the contrary, the harmonic balance method uses Fourier series, and the convergence rate is of exponential order. %Moreover, it is possible to notice that in a neighborhood of a Hopf bifurcation one harmonic is sufficient in order to get the best approximation of the periodic solution.
In our case, in general for the given parameters, about 50 harmonics are enough to get the desired accuracy which gives rise to similar numerical results for both methods.
% We can optimize the calculations, we control the following number of harmonic the excitation solutions, where the excitation is high, $I \in []$ using $50$ harmonic, after we can reduce the number of harmonics to $30$ for the interval for $I \in []$, for this case we have $3$ solutions see Fig. \ref{fig1}, for both solution have small amplitudes, we may use 30 harmonics, against by we keep $50$ harmonics for solution has large amplitude. For the solutions far  away Hopf bifurcation point $I \in []$, using 10 harmonics, then to the his neighborhood, a single harmonic is sufficient.
% %dans le reste on garde la methode de balance harmonique, comme etant un algorithme simple est efficace pour detecter les solution stable et instable et on trace les branches des points stationnaires et les solution periodique pour avoir le comportement general de systeme hh.

%Exploiting the numerical techniques introduced in the previous Section, we have derived the bifurcation diagram in Fig.1.

\begin{figure}[t!]
\begin{center}
\includegraphics[width=10cm,height=5cm]{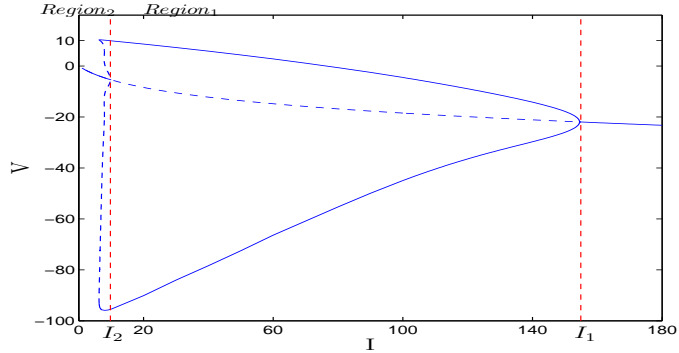}
\captionof{figure}{Branches of stable (solid line) and unstable (dotted line) periodic solutions of HH model. For each periodic solution the minimum and the maximum values of the potential $V$ over one period are represented. Depending on the values of $I$, two regions with different dynamical behaviors can be identified.} \label{fig1}
\end{center}
\end{figure}

In Fig. \ref{fig1} the bifurcation diagram for the HH model has been obtained by jointly and optimally exploiting the three methods presented above, that is shooting, collocation and HB methods. The stability analysis of the detected limit cycles has been carried out by the calculation of the Floquet multipliers, by applying the numerical algorithm proposed in \citep{farkas} to the approximated solution.

\begin{figure}[t!]
\begin{center}
\includegraphics[width=10cm,height=7cm]{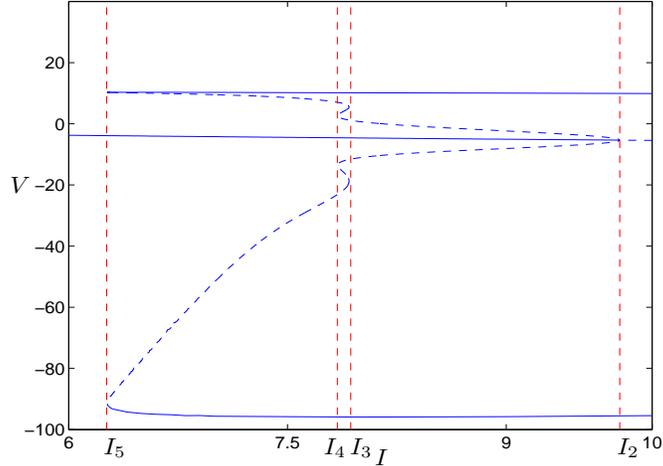}%{branche_2d_unstable.png}
\captionof{figure}{Zoom of region 2 of Fig \ref{fig1}. HH model exhibits one equilibrium point, one stable limit cycle (solid line) and up to 3 unstable ones (dotted lines).} \label{fig2}
\end{center}
\end{figure}

It is possible to see that the dynamical behavior of HH system can be decomposed in two main phases, depending on the value of the external current $I$. For $I > I_{2}=9.73749234$, there is only one equilibrium point and one stable periodic solution, that disappears through a Hopf bifurcation for $I=I_1=154.500$. The second phase is for $I \in~[0,~I_{2}]$ and it is more interesting, since its dynamical behavior is more complex and rather less understood. A zoom of Region 2 can be found in Fig. \ref{fig2}. It is easy to see that in this second region the system undergoes three saddle-node of cycles bifurcations at $I_3=7.92198549$, $I_4=7.84654752$ and $I_5=6.26490316$. They correspond to the knees of the bifurcation diagram and they consist in the collision and disappearance of two periodic solutions. Moreover, the system exhibits a period-doubling bifurcation \citep{Rinzel} at $I_6=7.92197768$ that can be detected by the joint application of the harmonic balance  method and the Floquet analysis.

\begin{figure}[!t]
    \begin{minipage}[l]{5cm}
        \centering
\includegraphics[scale=0.4]{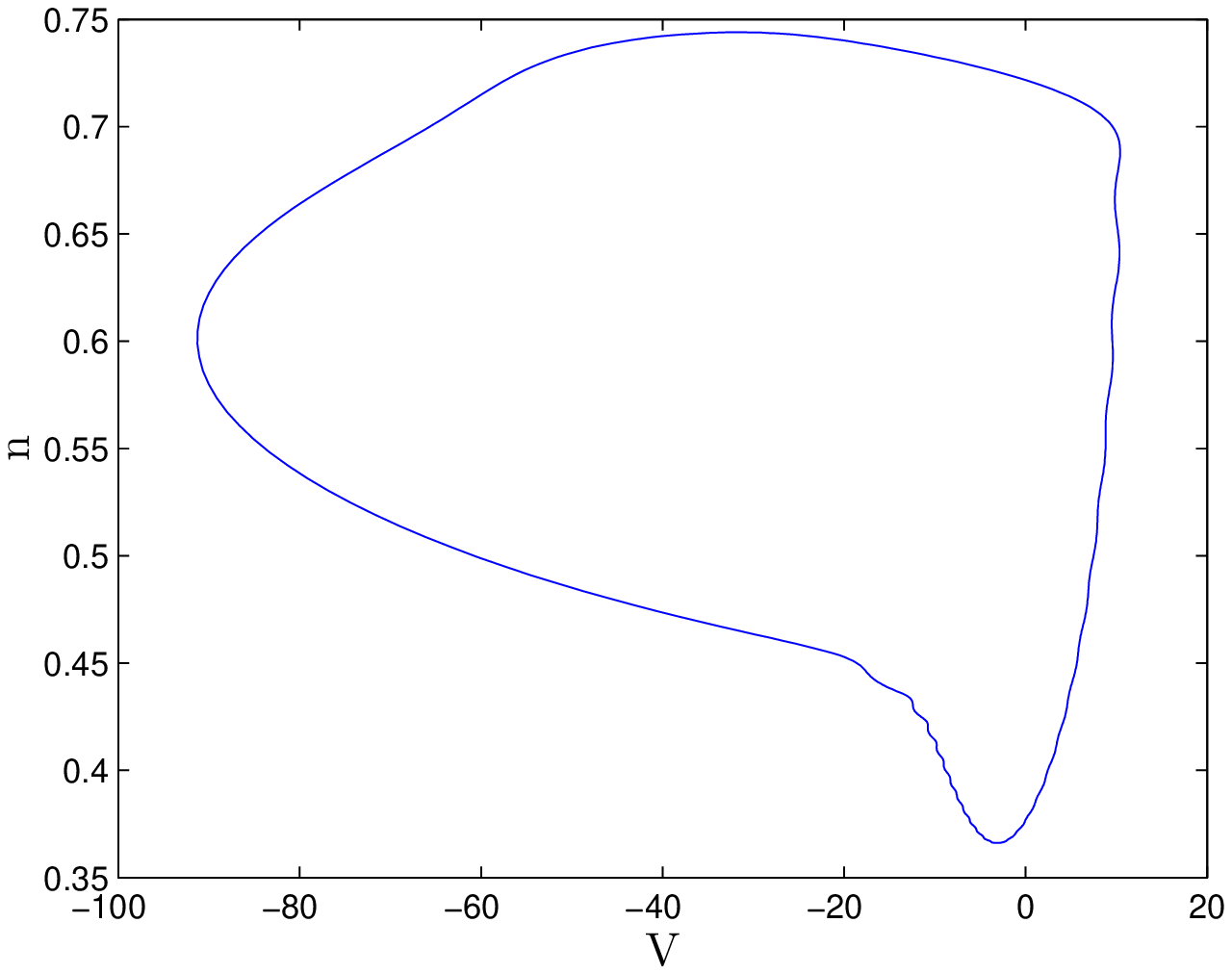}\\%{branche_2d_unstable.png}
                (a)
    \end{minipage}
\qquad
           \begin{minipage}[r]{5cm}
        \centering
        \includegraphics[scale=0.4]{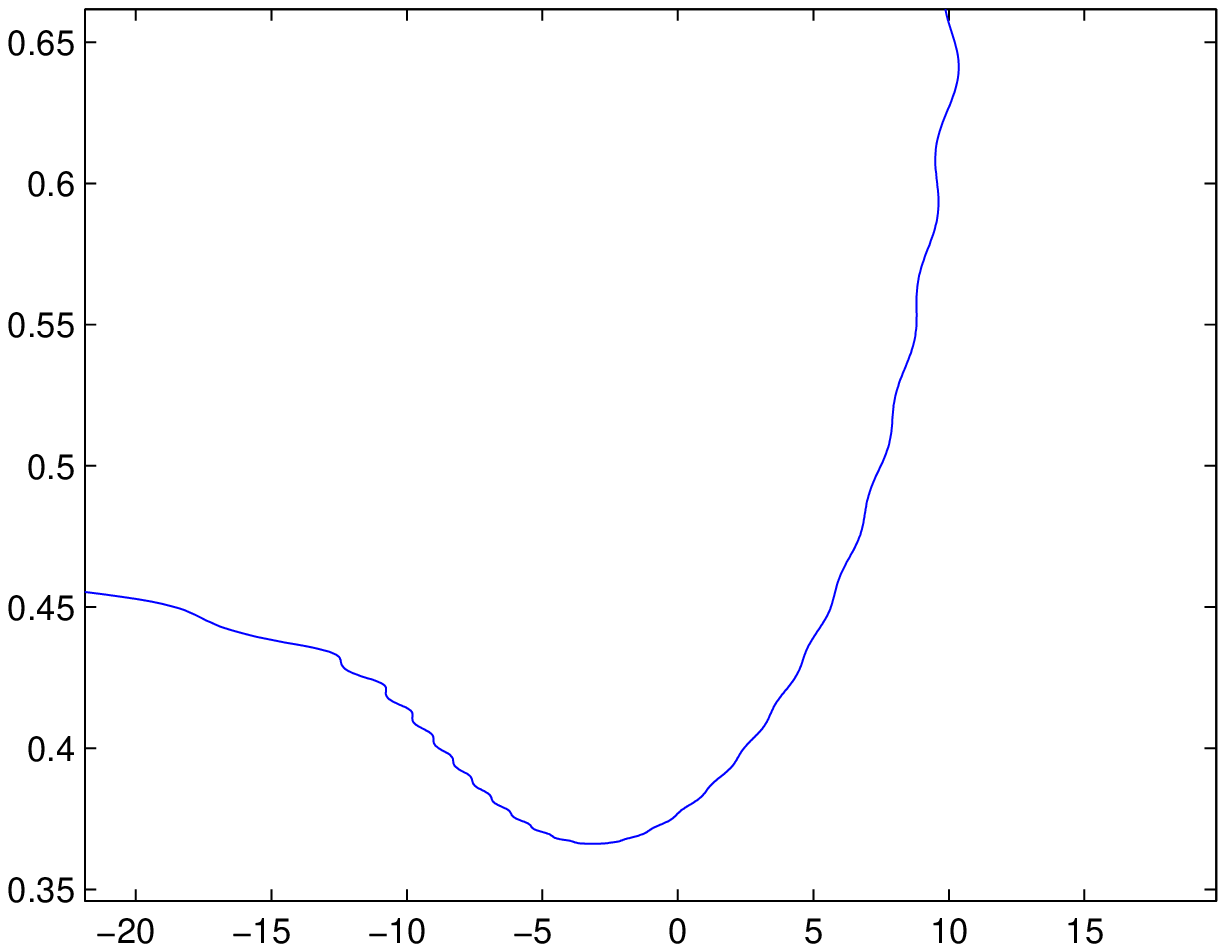}%{branche_2d_unstable.png}
\\
                (b)
    \end{minipage}
    \caption{(a) The stable periodic solution detected by the HB method for $I= 6.25$ exhibits the Gibbs phenomenon. (b) Zoom showing the small oscillations, sign of a non accurate approximation of the limit cycle, despite the exploitation of 50 harmonics.}
    \label{fig1bis}
\end{figure}

For $I$ close to $I_5$, the periodic solutions detected by the HB method exhibit the Gibbs phenomenon \citep{Urabe}, as it can be seen in Fig. \ref{fig1bis}, and this does not permit to accurately detect the saddle-node of  cycles bifurcation. Therefore only in this region shooting and collocation methods have been used in order to find the stable and unstable periodic solutions, respectively. The remaining of the diagram has been found via the HB method, by  choosing and controling the minimal number of required harmonics. In particular, in Region 2, 50 harmonics have been considered for the approximation of the high amplitude stable limit cycle, while for $I>7$ the unstable limit cycles required only 30 harmonics, and for $I > 8$ the number of harmonics can be gradually reduced since the unstable limit cycle becomes the more and more regular. Finally, close to the Hopf bifurcation at $I_2$ only one harmonics is sufficient to get the best approximation of the unstable periodic solution.

Therefore, we can conclude that HB method works very well in the region between $I_4$ and $I_2$, that is in the most interesting part of the diagram from a dynamical point of view. It permits to obtain the results in a more performant way with respect to the other methods.\\[0.1cm]

 In the following, we analyze more accurately those various bifurcations. \\[0.1cm]

\subsection{Analysis of the limit cycles bifurcations}
\begin{figure}[!t]
    \begin{minipage}[l]{5cm}
        \centering
\includegraphics[width=6cm,height=6cm]{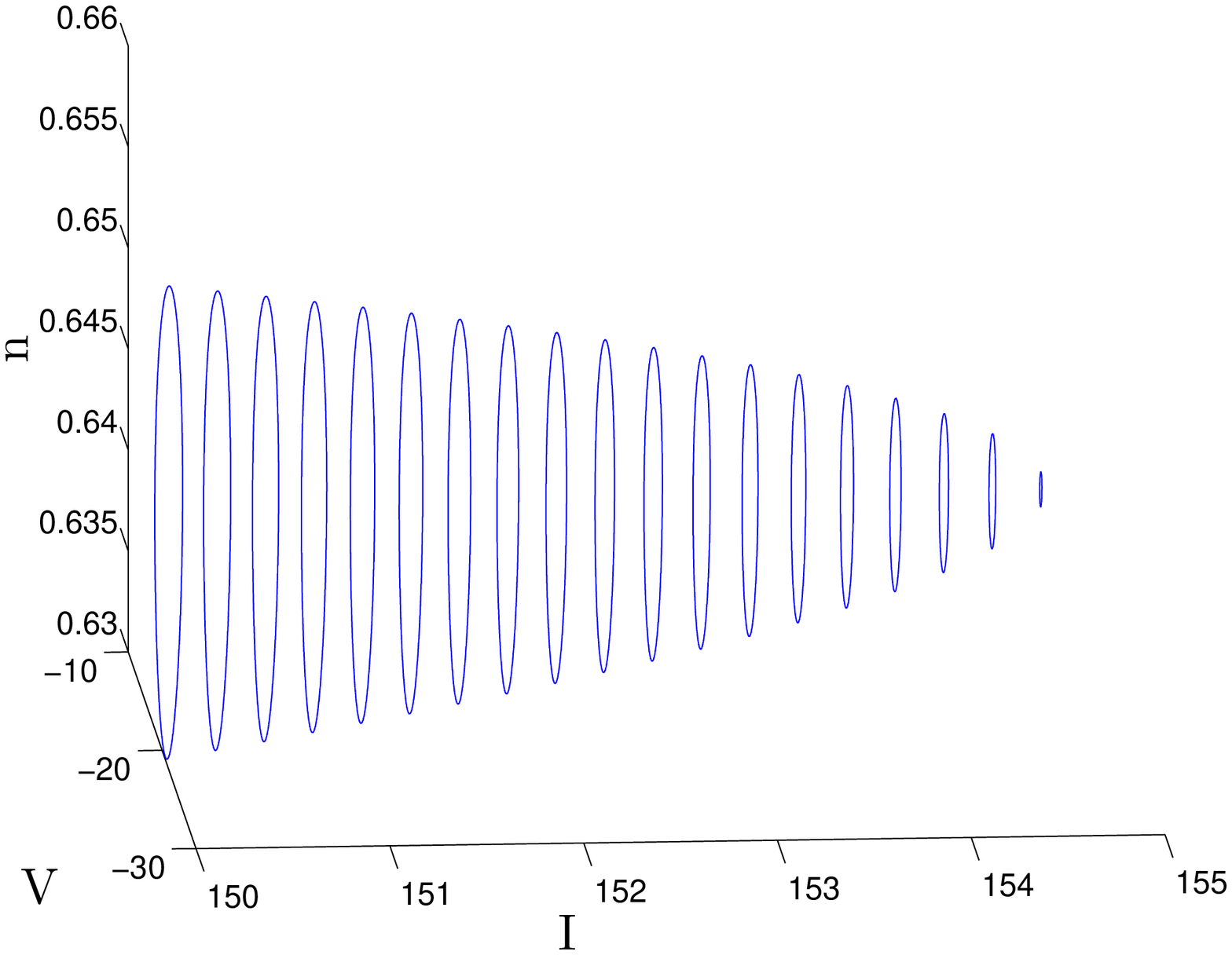}\\%{branche_2d_unstable.png}
                (a)
    \end{minipage}
           \begin{minipage}[r]{5cm}
        \centering
        \includegraphics[width=6cm,height=6cm]{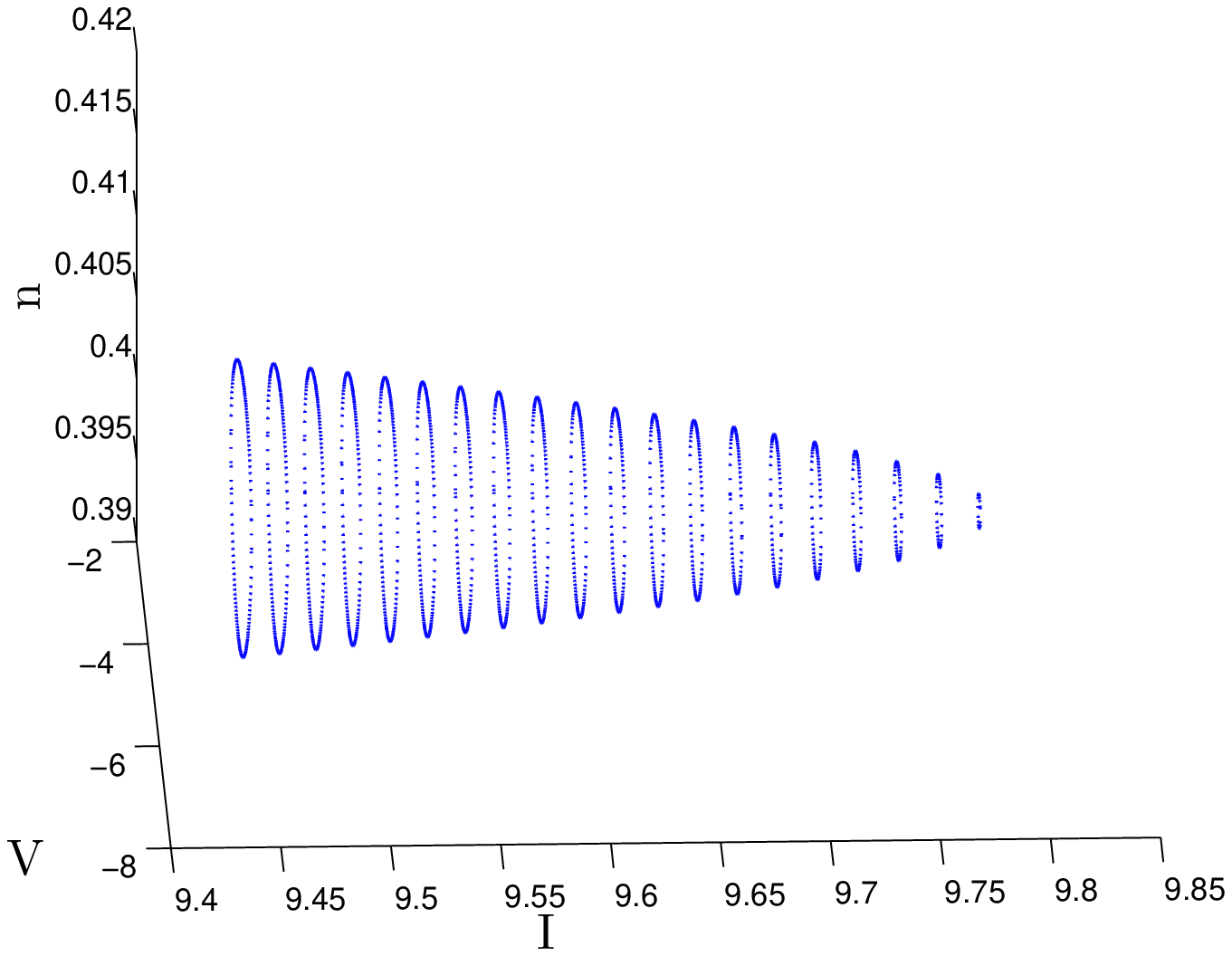}%{branche_2d_unstable.png}
\\
                (b)
    \end{minipage}
    \caption{Branches of (a) stables and (b) unstables periodic solutions for different values of $I$, in a neighborhood of (a) $I_1$ and (b) $I_2$, respectively.}
    \label{figHopf1}
\end{figure}

\textbf{Hopf bifurcations.} In this paragraph, we are interested in Hopf bifurcations, that take place at $I=I_{1}=154.500$ and $I_{2}=9.73749234$. A view into  $(V,n,I)$ space projection of stable and unstable periodic solutions, in a neighborhood of the two Hopf bifurcations, are shown in Fig. \ref{figHopf1}. These results suitably match with the theoretical results proved in \citep{Hassard,Chaos-HH}.

% \begin{figure}[h!]
% \begin{center}
% \includegraphics[width=6cm,height=6cm]{brache_bfh_2}%{branche_2d_unstable.png}
% \captionof{figure}{The branch of the stable periodic solution} \label{fig1}
% \end{center}
% \end{figure}
%
% \begin{figure}[h!]
% \begin{center}
% \includegraphics[width=6cm,height=6cm]{branchehopf_unsta}%{branche_2d_unstable.png}
% \captionof{figure}{Branch of unstable periodic solutions} \label{fig1}
% \end{center}
% \end{figure}
%

\begin{figure}[!t]
    \begin{minipage}[l]{5cm}
        \centering
\includegraphics[width=5cm,height=5cm]{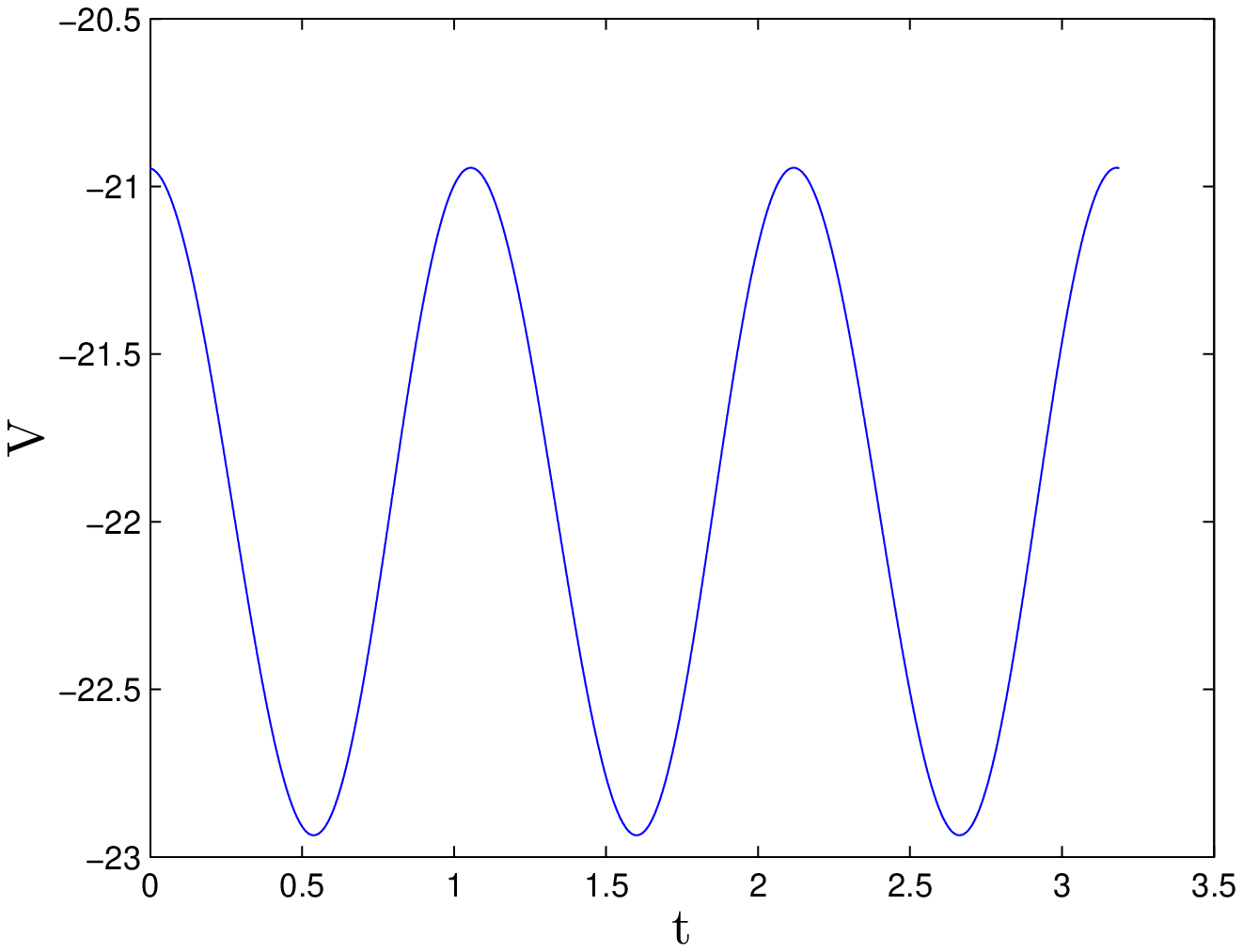}\\
                (a)
    \end{minipage}
           \begin{minipage}[r]{5cm}
        \centering
        \includegraphics[width=5cm,height=5cm]{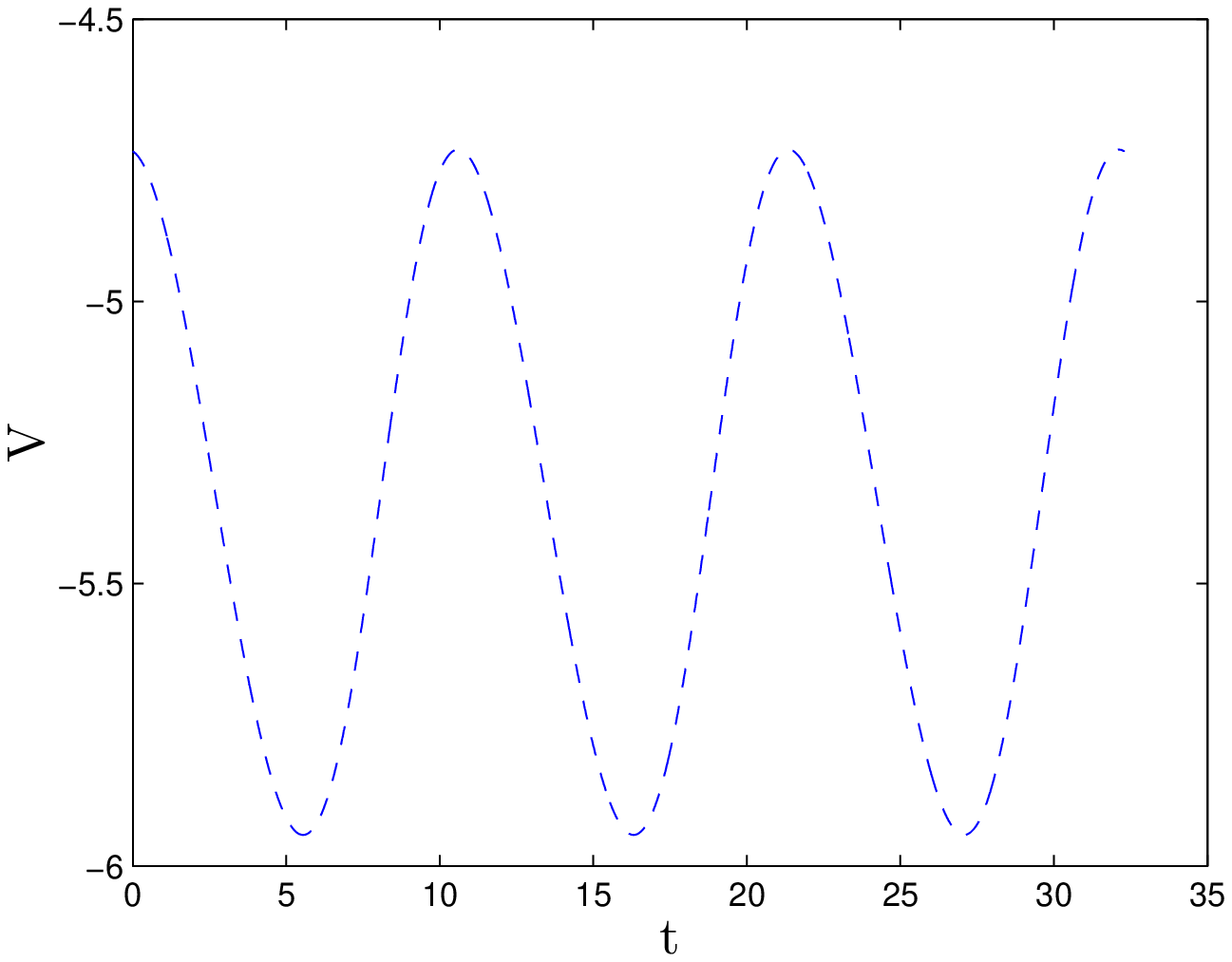}\\
                (b)
    \end{minipage}
    \caption{ (a) Stable periodic solution for $I=152.2500$ and (b) unstable periodic solution for $I=9.71889$.}
    \label{figHopf2}
\end{figure}

It is worth noting that in both cases over a large interval of $I$ close to the Hopf bifurcations, the periodic solutions are almost sinusoidal (see Fig. \ref{figHopf2}). Therefore, only one or two harmonics are needed to conveniently approximate this solution via the harmonic balance method. On the contrary, the collocation method still requires a huge number of nodes, so in this case the nonlinear system to solve is still of high dimension.\\[0.1cm]

% In our case the spectral method is independent of the instability of solutions, that is to say that either the bifurcation  of type sub critic $I=I_1$ or critical $I=I_2$ there needs to a one harmonic ( number of node equal 3)  to approximate the solution
% However the collocation method still requires (needs) an interesting number of nodes, even for the sine function case.
% The fig and fig explore the sine function near to Hopf bifurcation.

\begin{figure}[!t]
    \begin{minipage}[l]{5cm}
        \centering
\includegraphics[width=5cm,height=5cm]{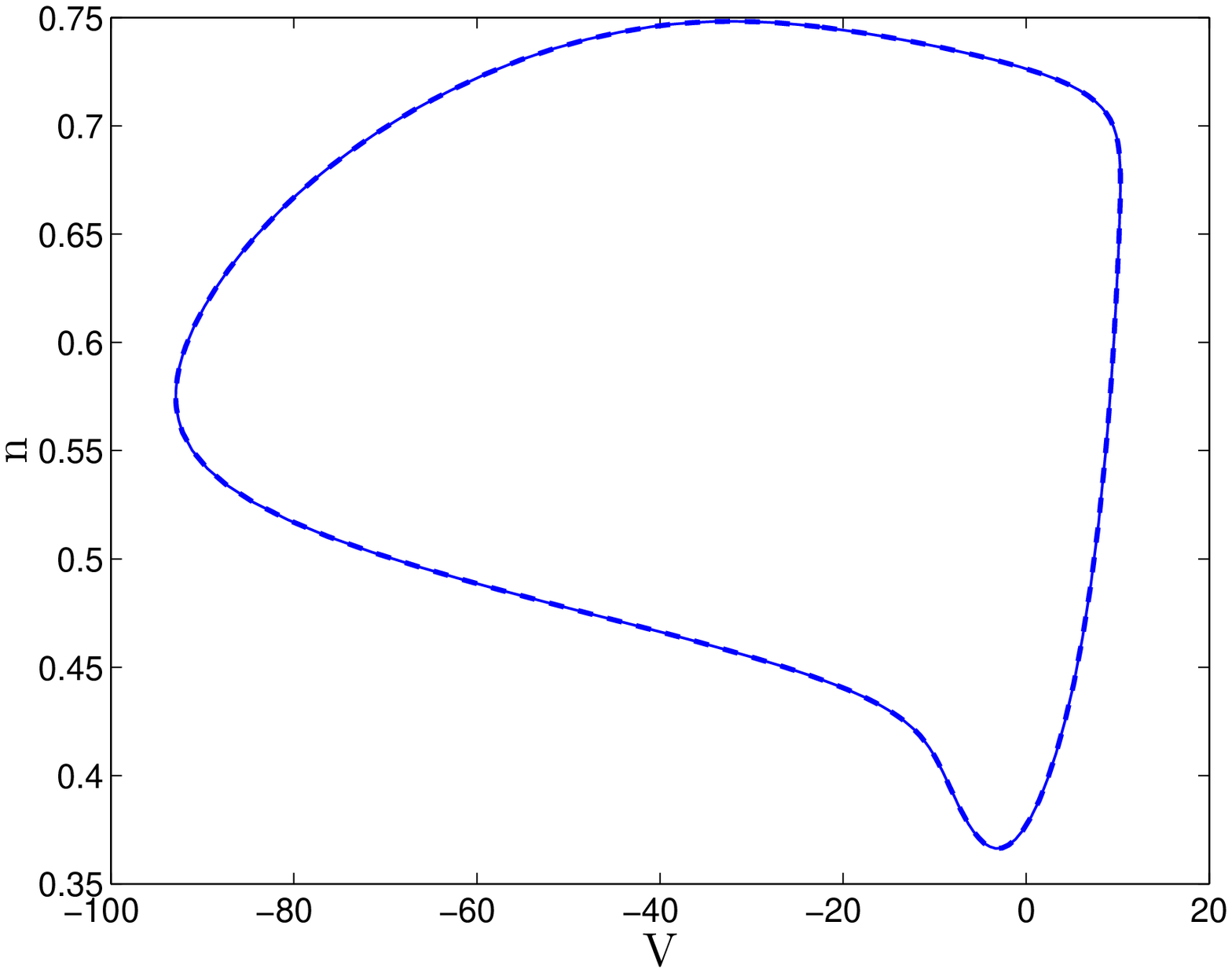}\\

                (a)
    \end{minipage}
\quad
           \begin{minipage}[r]{5cm}
        \centering
\includegraphics[width=5cm,height=5cm]{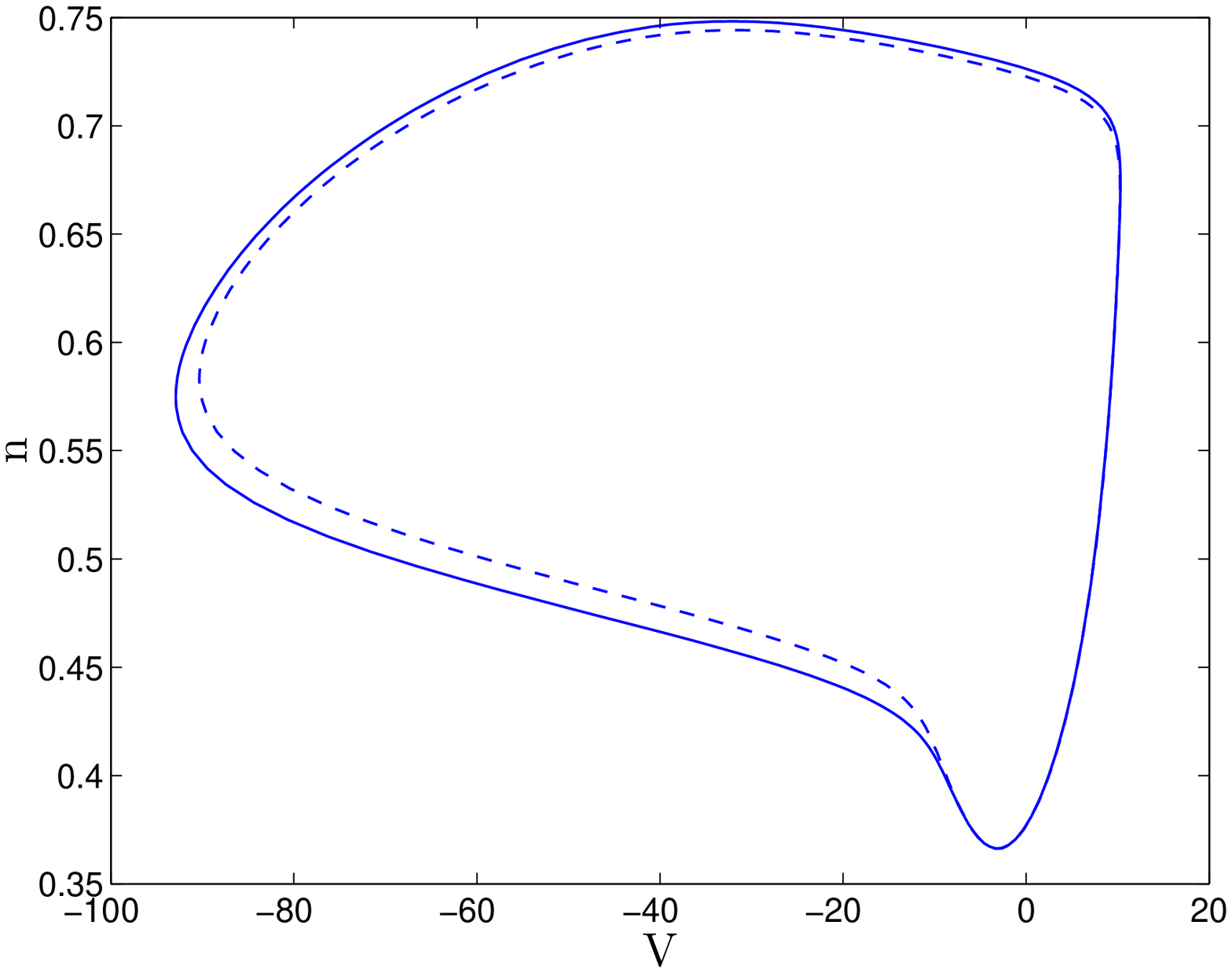}\\
                (b)
    \end{minipage}
    \caption{Stable (solid line) and unstable (dashed line) limit cycles near the first saddle-node of cycles bifurcation, for (a) $I= 6.2649$ both solutions are almost coincidents, and for (b) $I= 6.2716$.}
\label{figSN1}\end{figure}

\begin{figure}[!t]
    \begin{minipage}[l]{5cm}
        \centering
\includegraphics[width=5cm,height=5cm]{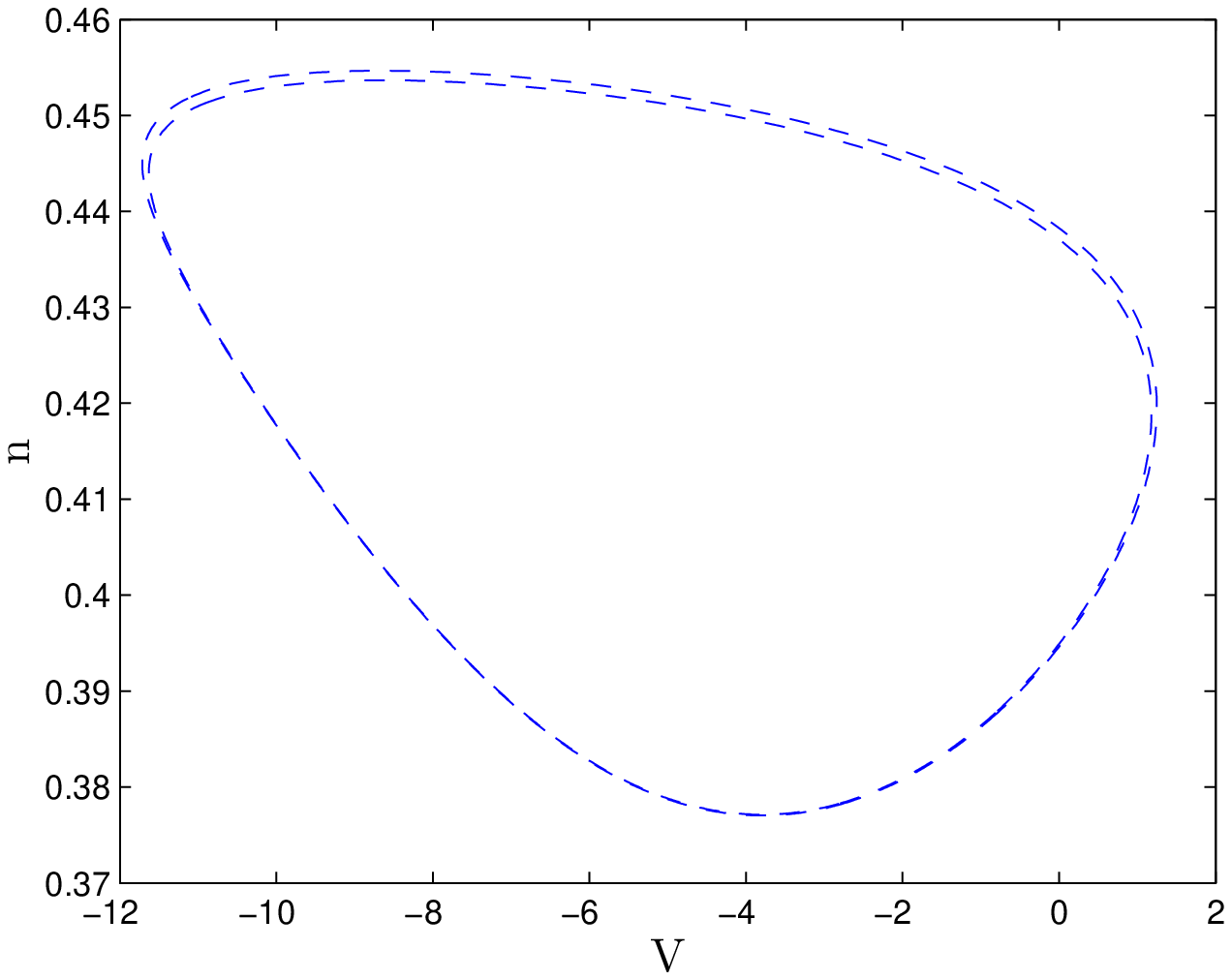}\\
                (a)
    \end{minipage}
\quad
           \begin{minipage}[r]{5cm}
        \centering
        \includegraphics[width=5cm,height=5cm]{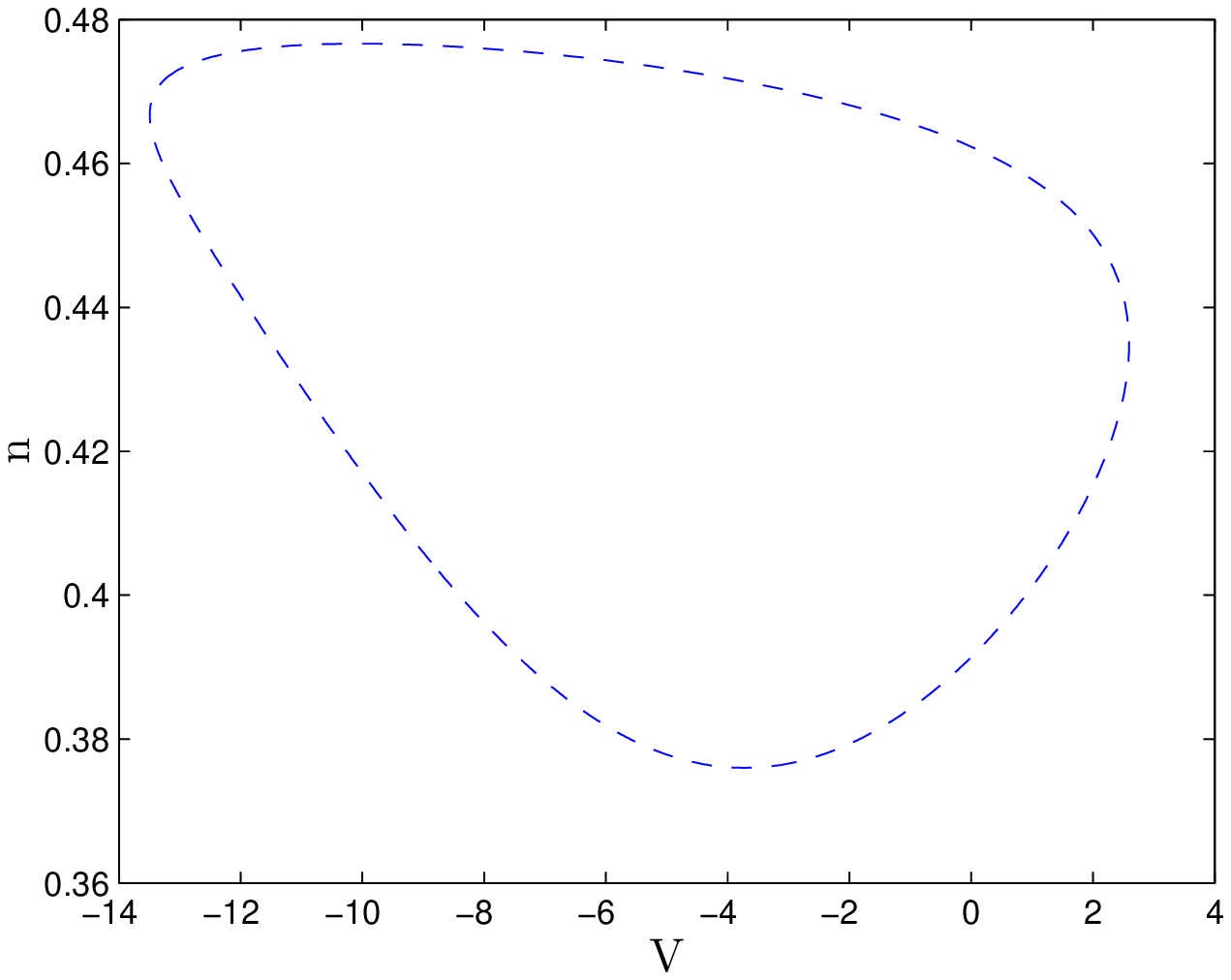}%{branche_2d_unstable.png}
\\
                (b)
    \end{minipage}
    \caption{Projection of two unstables limit cycles on the $(V,n)$ plane for (a) $I=7.92198548\lesssim I_3$ and (b)   $I=I_3=7.92198549$.  }
\label{figSN2}\end{figure}

\begin{figure}[!t]
    \begin{minipage}[l]{5cm}
        \centering
\includegraphics[width=5cm,height=5cm]{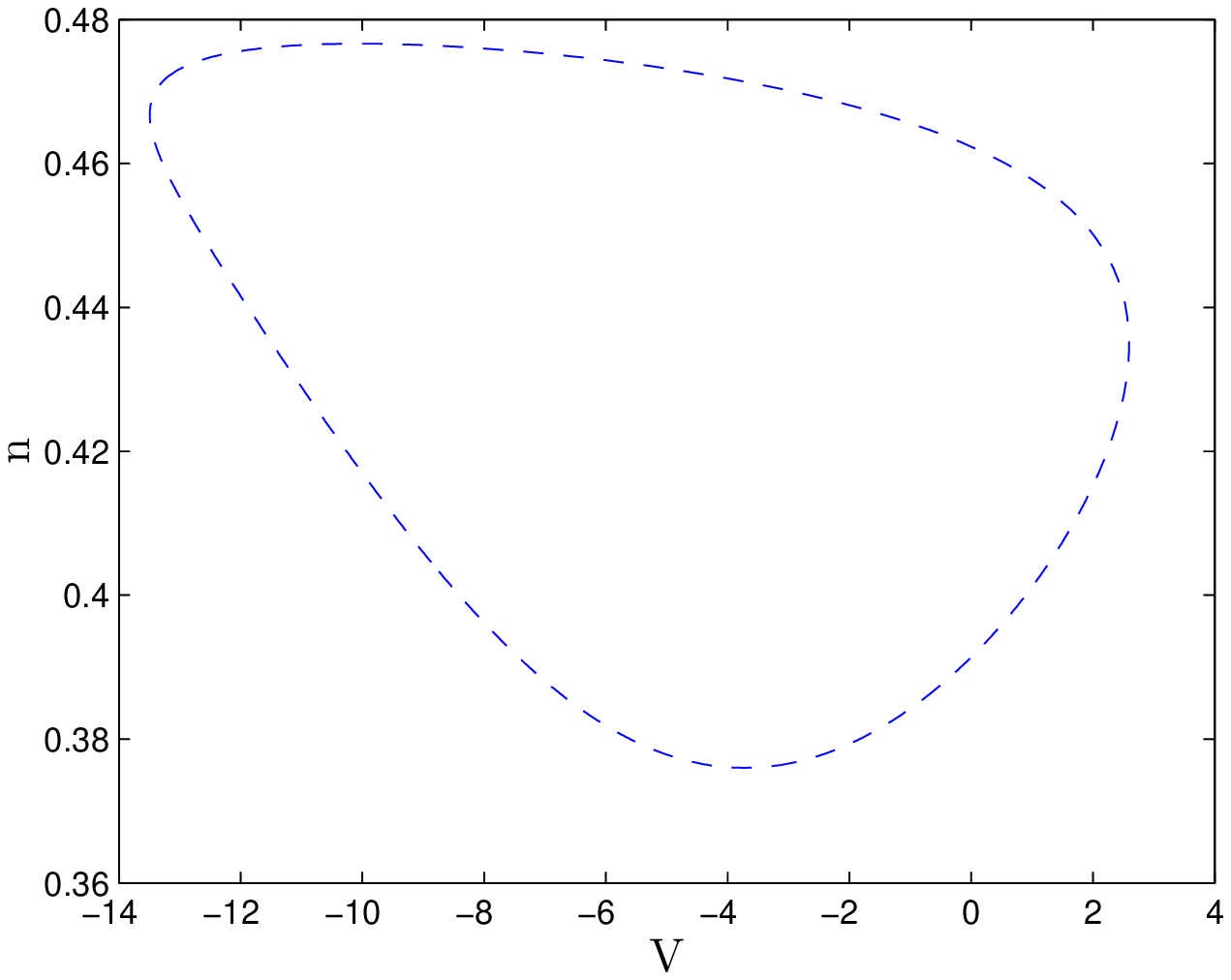}\\%{branche_2d_unstable.png}
                (a)
    \end{minipage}
           \begin{minipage}[r]{5cm}
        \centering
        \includegraphics[width=5cm,height=5cm]{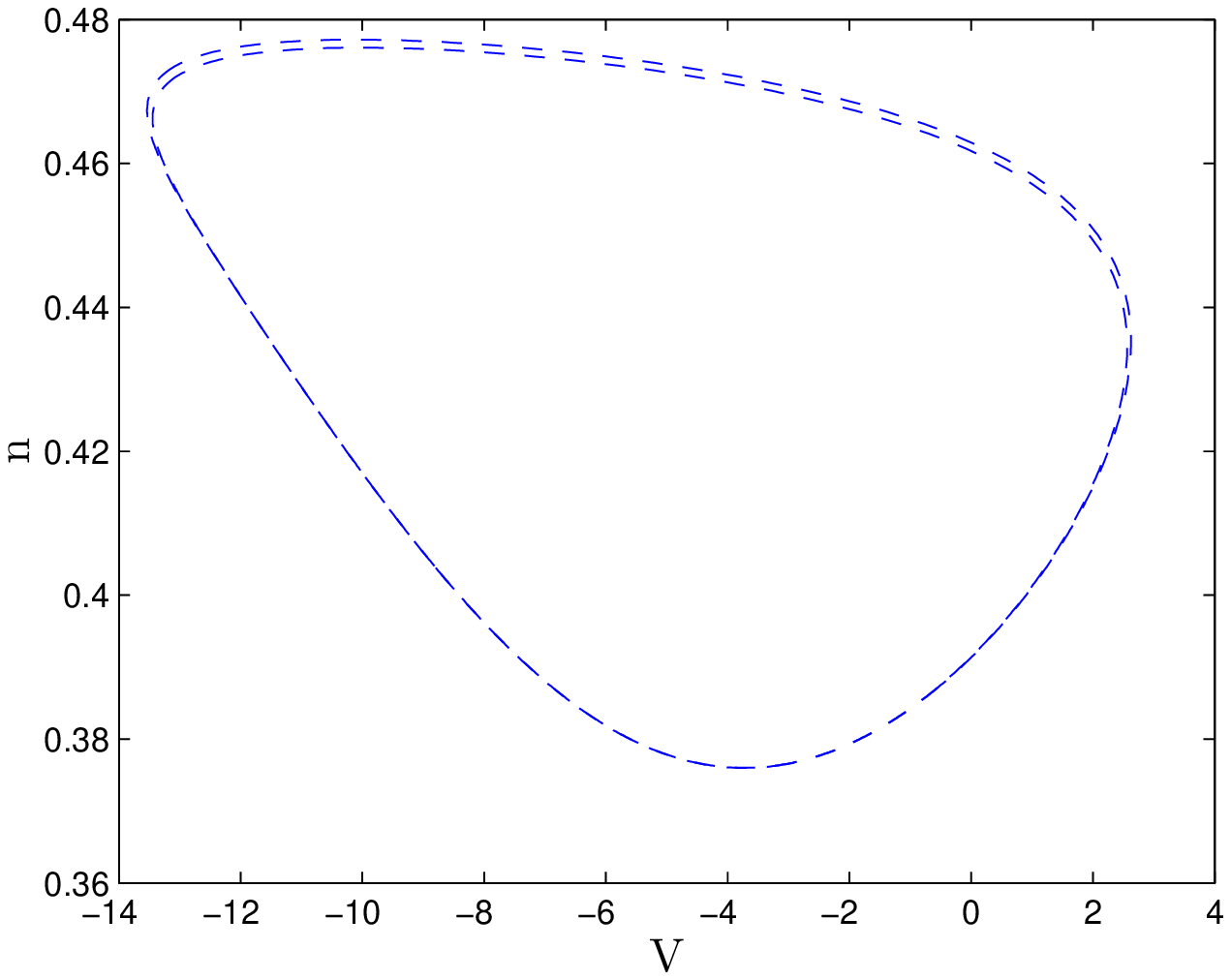}%{branche_2d_unstable.png}
\\
                (b)
    \end{minipage}
    \caption{Projection of two unstable limit cycles on the $(V,n)$ plane for (a) $I=I_4=7.84654752$ and (b) $I=7.84654876\lesssim I_4$.  }
\label{figSN3}\end{figure}

\noindent\textbf{Saddle node of cycles bifurcations}
In our case, there are two types of saddle node of cycles bifurcation: for $I = I_5=6.26490316$ we have a simultaneous appearence of two limit cycles (one stable and the other unstable), while at $I = I_3=7.92198549$ and $I=I_4=7.84654752$ we have the collision of two unstable periodic solutions (see Figs. \ref{figSN1}, \ref{figSN2} and \ref{figSN3}). For detecting such bifurcations, we use the Floquet analysis, by searching when an additional Floquet multiplier crosses the unit circle in $+1$.

% Rinzel use other technique based to $\frac{dI}{dw}=0$, where $w=\frac{2 \pi}{T}$, it is the frequency.

\begin{figure}[p!]
   \begin{minipage}[l]{4cm}
     \centering \includegraphics[width=4.5cm,height=4.5cm,scale=0.35]{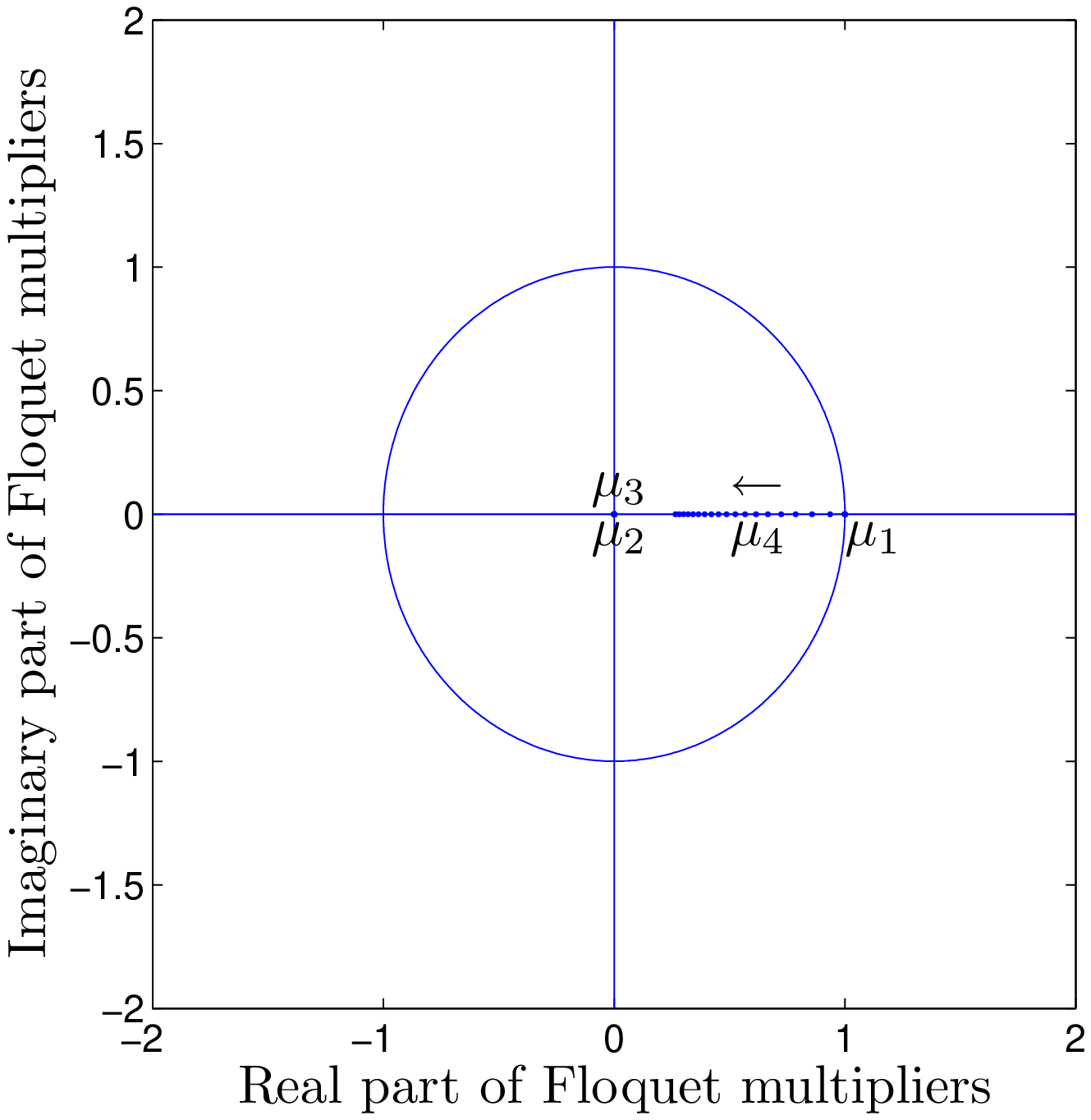}\\
      %\hspace*{3.5cm}
        (a)\vspace{3ex}
    \end{minipage}\hfill
   \begin{minipage}[r]{4cm}  
      \centering \includegraphics[width=4.5cm,height=4.5cm,scale=0.35]{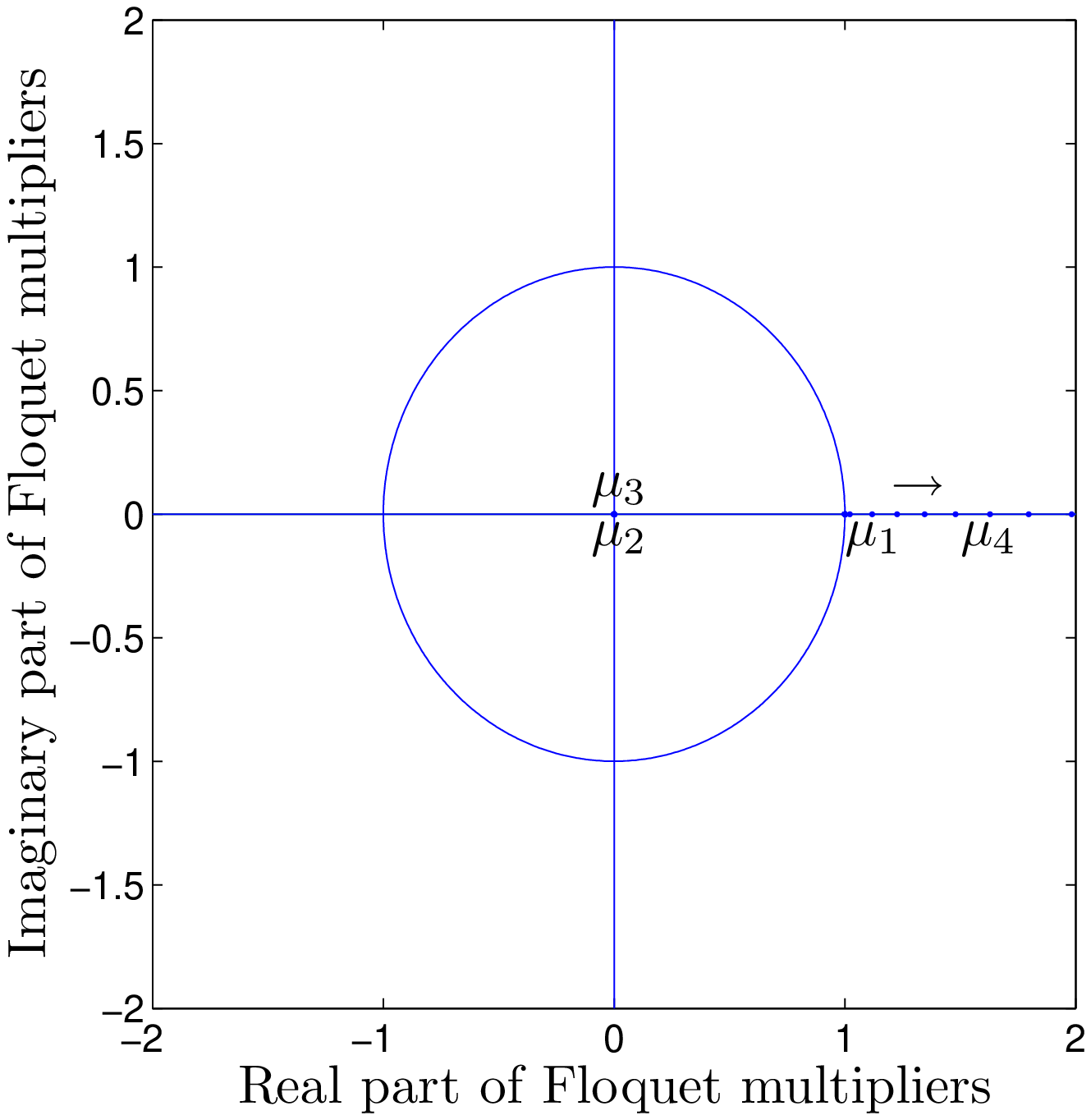}
      %\hspace*{3.5cm}
        (b)\vspace{3ex}
   \end{minipage}\\[-0.4cm]
   \begin{minipage}[l]{4cm}
      \centering \includegraphics[width=4.5cm,height=4.5cm,scale=0.35]{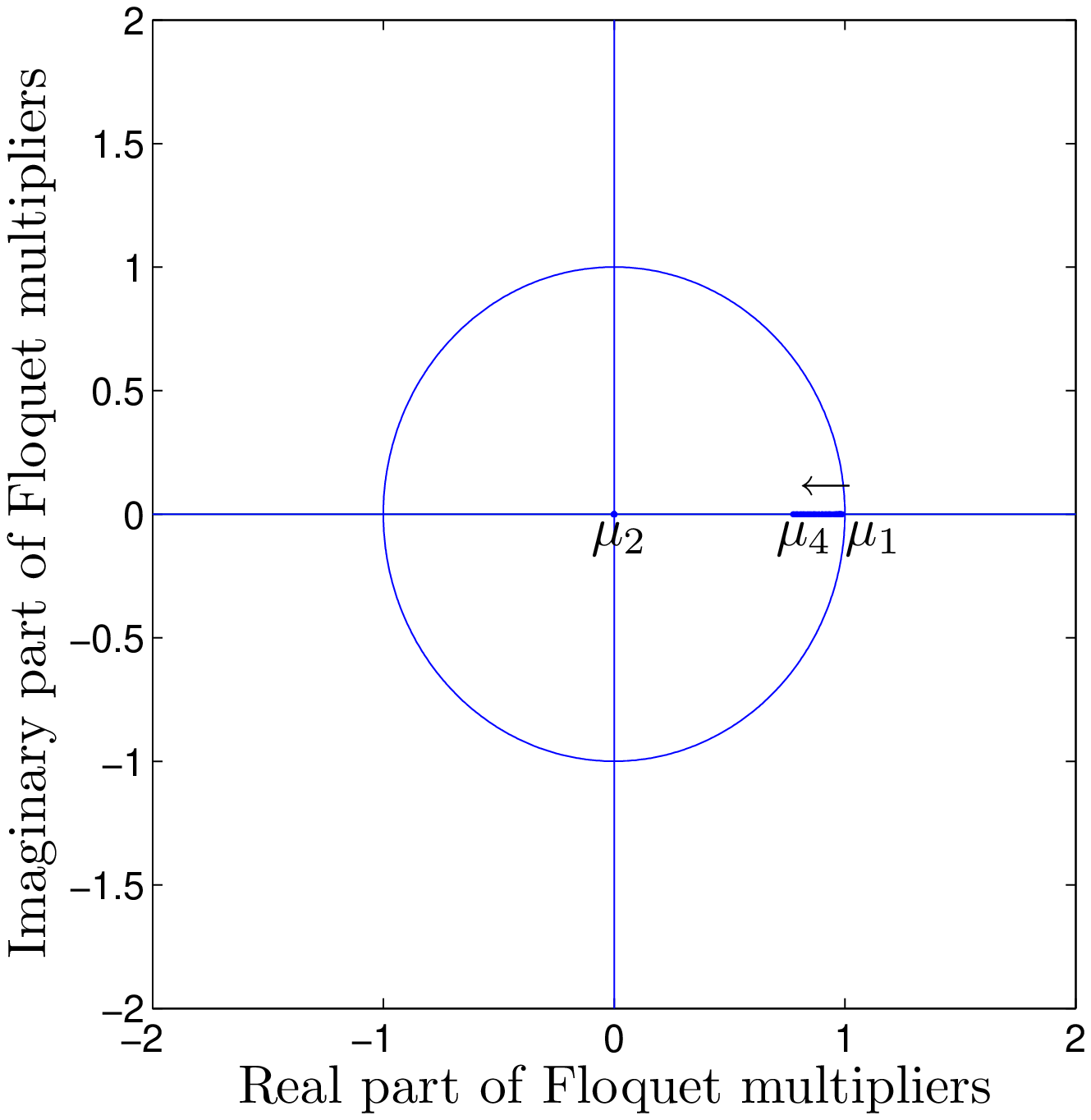}
       %\hspace*{3.5cm}
       (c)\vspace{3ex}
   \end{minipage}\hfill
   \begin{minipage}[r]{4cm}  
      \centering \includegraphics[width=4.5cm,height=4.5cm,scale=0.35]{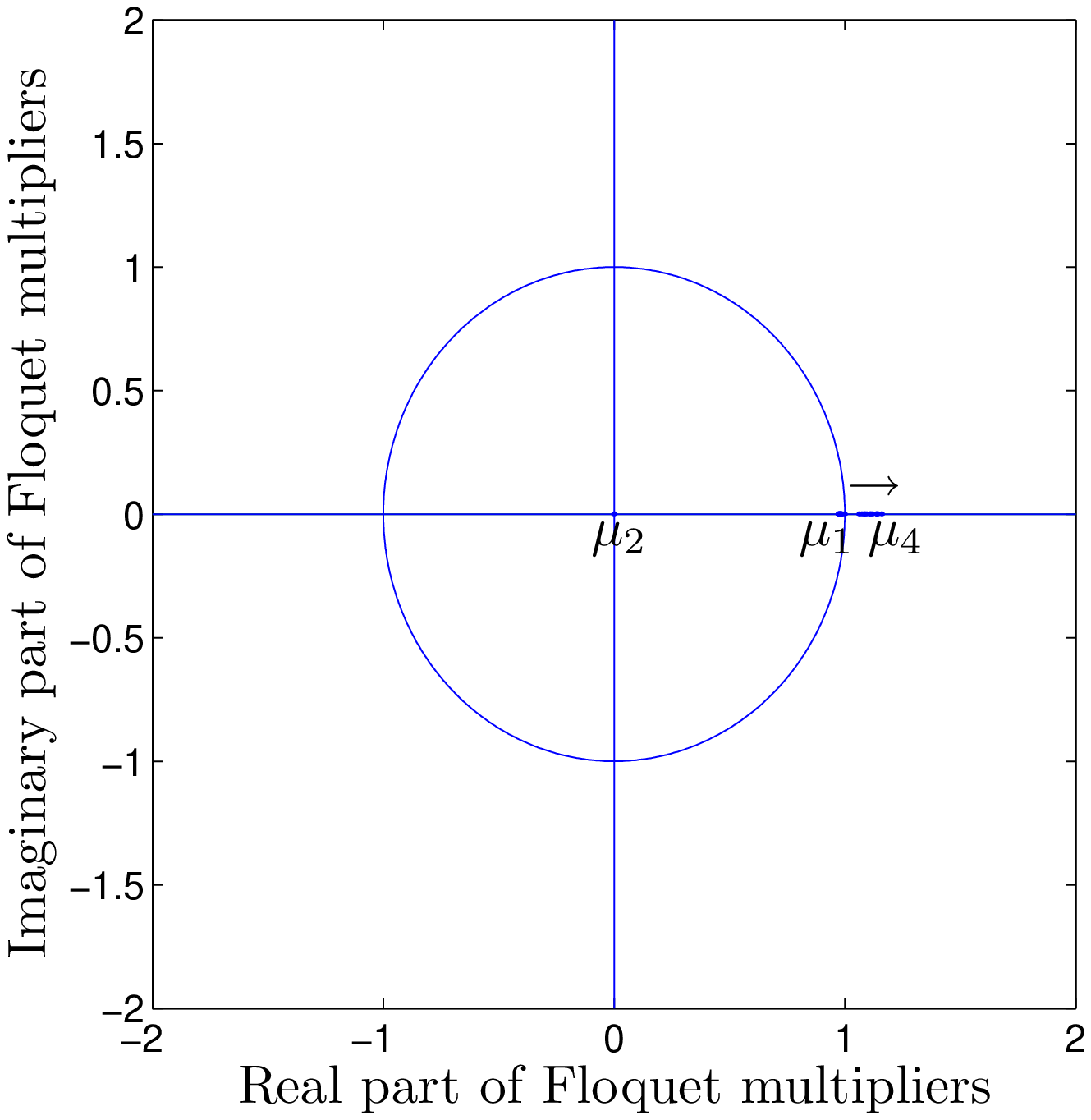}
      %  \hspace*{3.5cm}
           (d)\vspace{3ex}
   \end{minipage}\\[-0.4cm]
   \begin{minipage}[l]{4cm}
      \centering \includegraphics[width=4.5cm,height=4.5cm,scale=0.35]{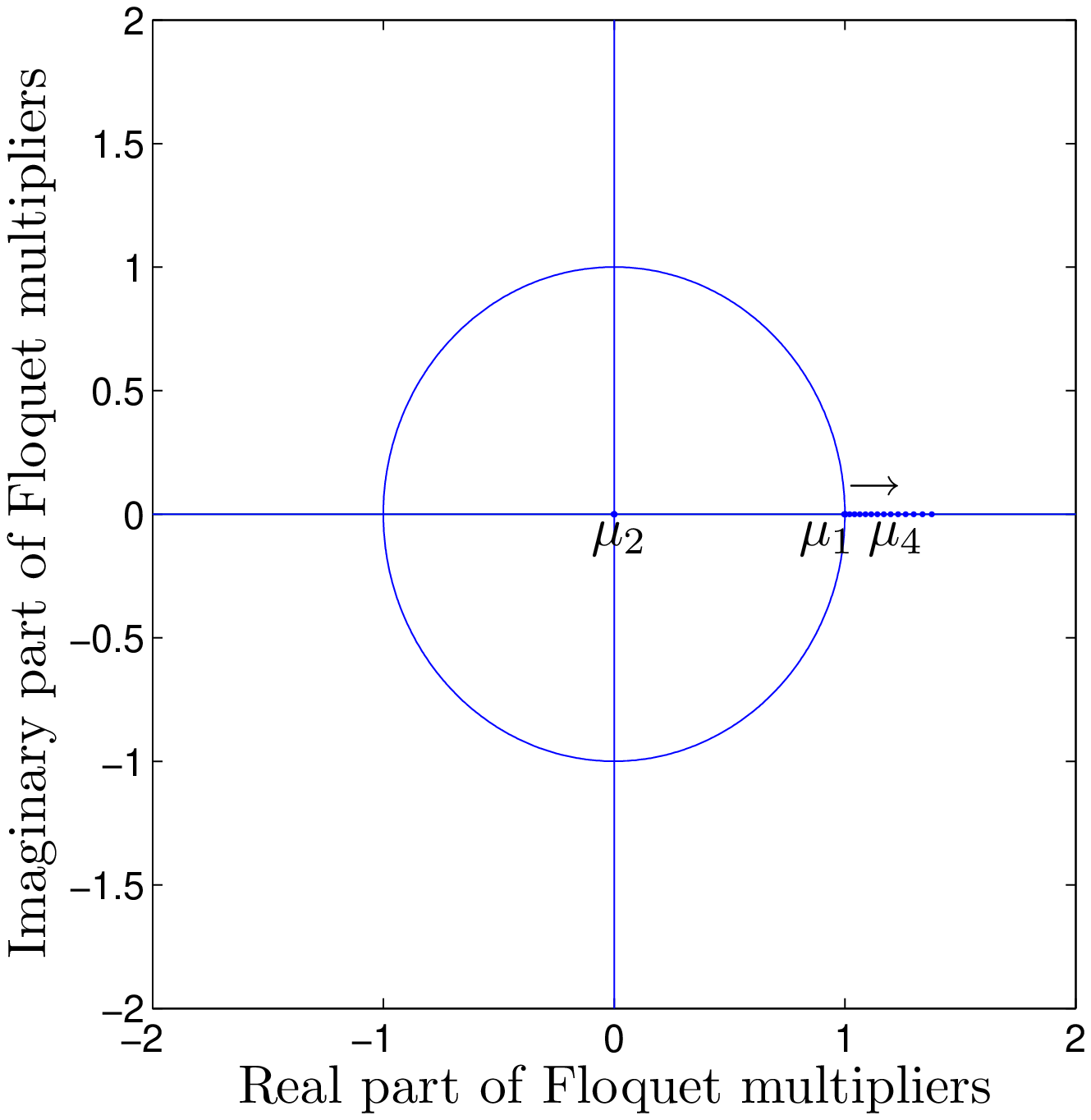}
     % \hspace*{3.5cm}
           (e)
   \end{minipage}\hfill
   \begin{minipage}[r]{4cm}  
      \centering \includegraphics[width=4.5cm,height=4.5cm,scale=0.35]{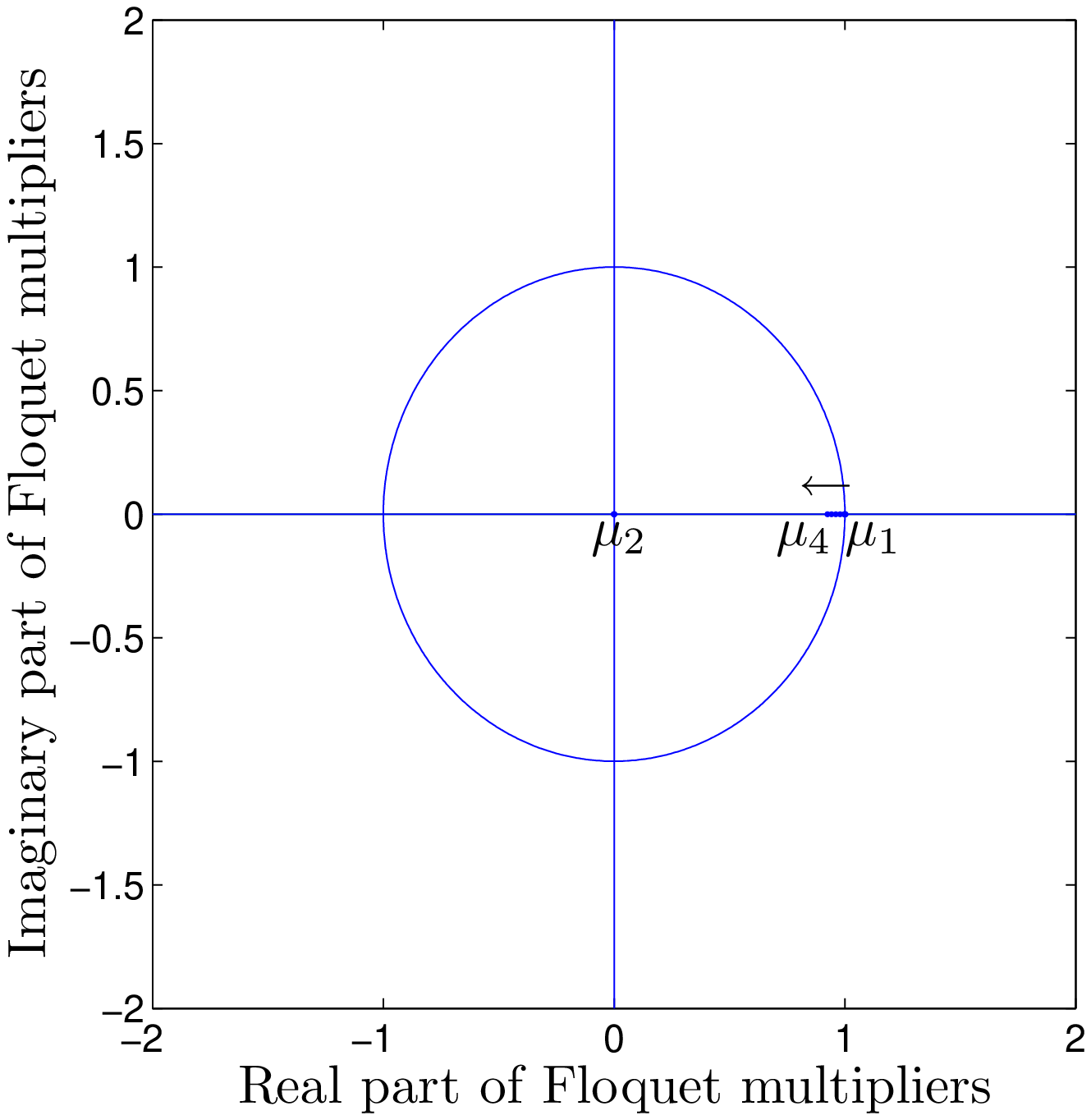}
       %\hspace*{3.5cm}
          (f)
   \end{minipage}

% \begin{figure}[p!]
%    \begin{minipage}[l]{3.5cm}
%      \centering \includegraphics[bb=137 208 540 622,scale=0.35]{mul_1bb.eps}\\
%       %\hspace*{3.5cm}
%         (a)\vspace{4ex}
%     \end{minipage}\hfill
%    \begin{minipage}[r]{3.5cm}  
%       \centering \includegraphics[bb=137 208 540 622,scale=0.35]{mul_2bb.eps}
%       %\hspace*{3.5cm}
%         (b)\vspace{4ex}
%    \end{minipage}\\[-0.4cm]
%    \begin{minipage}[l]{3.5cm}
%       \centering \includegraphics[bb=137 208 540 622,scale=0.35]{mul_unst1bb.eps}
%        %\hspace*{3.5cm}
%        (c)\vspace{4ex}
%    \end{minipage}\hfill
%    \begin{minipage}[r]{3.5cm}  
%       \centering \includegraphics[bb=137 208 540 622,scale=0.35]{mul_unst2bb.eps}
%       %  \hspace*{3.5cm}
%            (d)\vspace{4ex}
%    \end{minipage}\\[-0.4cm]
%    \begin{minipage}[l]{3.5cm}
%       \centering \includegraphics[bb=137 208 540 622,scale=0.35]{mul_unst3bb.eps}
%      % \hspace*{3.5cm}
%            (e)
%    \end{minipage}\hfill
%    \begin{minipage}[r]{3.5cm}  
%       \centering \includegraphics[bb=137 208 540 622,scale=0.35]{mul_unst4bb.eps}
%        %\hspace*{3.5cm}
%           (f)
%    \end{minipage}
%    
  \caption{ (a)-(b) Floquet multipliers for the stable limit cycles and unstable limit cycles, respectively, associated to the first saddle node of cycles bifurcation for $I\in [6.2792,6.7872]$. As $I$ increases, in (a) the multiplier $\mu_4$ starts from the value +1 and then enters in the unit circle, while in (b) the multiplier $\mu_4$, starts to the value +1 and becomes bigger and bigger. (c)-(d) Floquet multipliers for the two unstable limit cycles associated to the second saddle node of cycles bifurcation for $I\in [7.921985465,7.921985491]$. 
Here, in both cases, the third multiplier is outside the unit circle (this makes the limit cycle unstable) and is not shown, since it takes very high values with respect to the others. As in the previous case, as $I$ decreases, the multiplier $\mu_4$ starts from the value +1 and either (c) enters in the unit circle, or (d) takes higher and higher values.  (e)-(f) Floquet multipliers for the two unstable limit cycles associated to the third saddle node of cycles bifurcation for $I\in [7.846557778, 7.846616827]$. Also in this case, for both limit cycles, the  third multiplier is not represented. As $I$ increases, the multiplier $\mu_4$ starts from the value 
+1 and either (e) escapes from, or (f) enters in the unit circle.}\label{figSNmult}
   
\end{figure}

The Floquet multipliers for these three cases are represented in Fig.\ref{figSNmult}. It is possible to see that a multiplier leaves or enters in the unit circle through $+1$.

% \begin{figure}[p!]
%    \begin{minipage}[l]{5cm}
%      \centering \includegraphics[width=5.5cm,height=5.5cm,scale=0.35]{mul_1bb.eps}\\
%       %\hspace*{3.5cm}
%         (a)\vspace{4ex}
%     \end{minipage}\hfill
%    \begin{minipage}[r]{5cm}  
%       \centering \includegraphics[width=5.5cm,height=5.5cm,scale=0.35]{mul_2bb.eps}
%       %\hspace*{3.5cm}
%         (b)\vspace{4ex}
%    \end{minipage}
%    \begin{minipage}[l]{5cm}
%       \centering \includegraphics[width=5.5cm,height=5.5cm,scale=0.35]{mul_unst1bb.eps}
%        %\hspace*{3.5cm}
%        (c)\vspace{4ex}
%    \end{minipage}\hfill
%    \begin{minipage}[r]{5cm}  
%       \centering \includegraphics[width=5.5cm,height=5.5cm,scale=0.35]{mul_unst2bb.eps}
%       %  \hspace*{3.5cm}
%            (d)\vspace{4ex}
%    \end{minipage}
%    \begin{minipage}[l]{5cm}
%       \centering \includegraphics[width=5.5cm,height=5.5cm,scale=0.35]{mul_unst3bb.eps}
%      % \hspace*{3.5cm}
%            (e)
%    \end{minipage}\hfill
%    \begin{minipage}[r]{5cm}  
%       \centering \includegraphics[width=5.5cm,height=5.5cm,scale=0.35]{mul_unst4bb.eps}
%        %\hspace*{3.5cm}
%           (f)
%    \end{minipage}
 
\begin{figure}[t!]
\begin{center}
\includegraphics[scale=0.4]{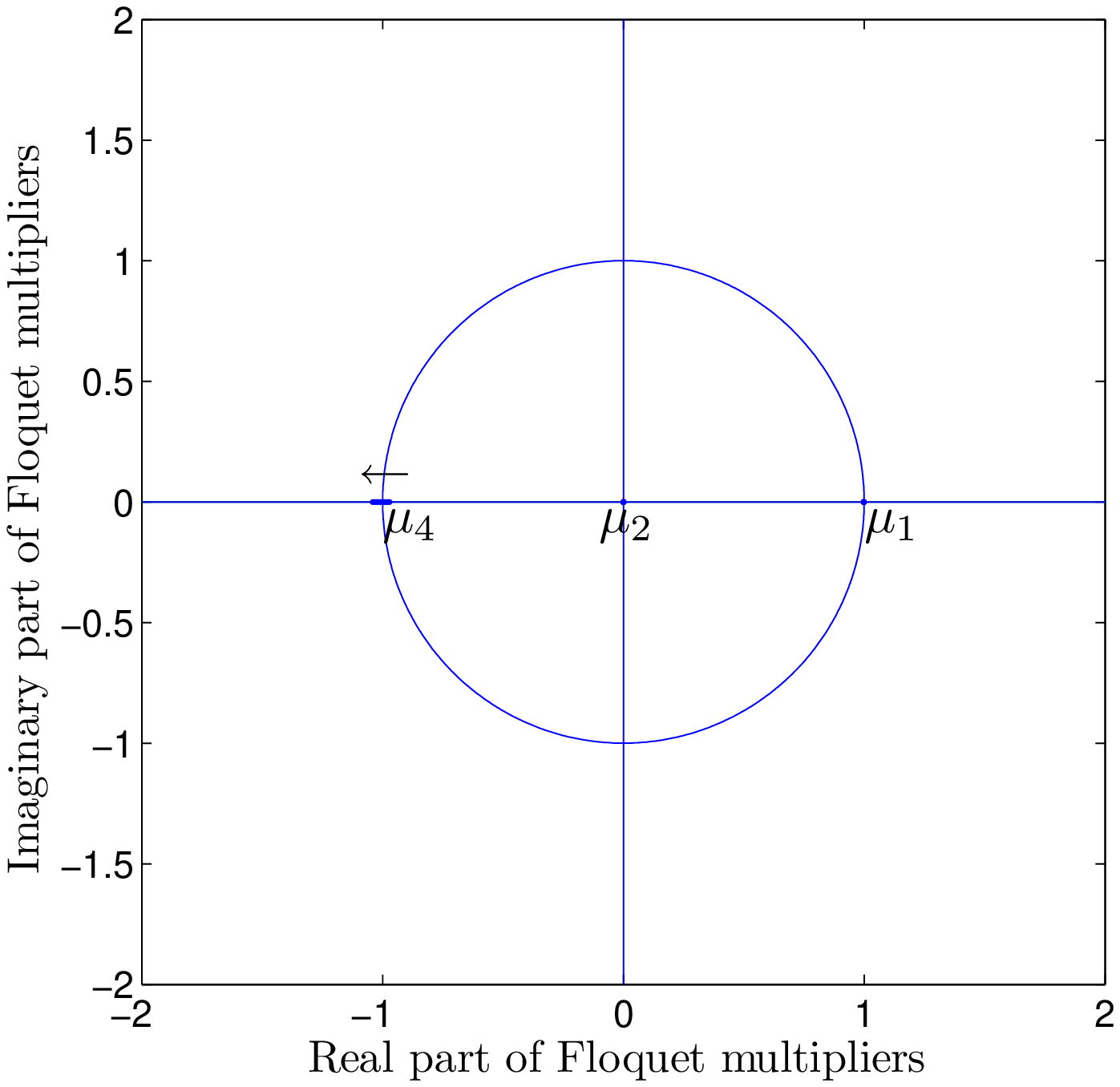}%{branche_2d_unstable.png}
\captionof{figure}{Floquet multipliers near the period-doubling bifurcation for different values of $I\in [ 7.92197743, 7.92197799 ]$. By decreasing $I$, the multiplier $\mu_4$ crosses the unit cycle through $-1$.} \label{figPD}
\end{center}
\end{figure}

\vspace{0.2cm}

\textbf{Period doubling bifurcation} Finally, in this paragraph, we consider the period-doubling bifurcation. By exploiting the Floquet analysis, we can easily detect this bifurcation since in this case a Floquet multiplier crosses the unit circle through $-1$, as it is shown in Fig. \ref{figPD}. Table \ref{table} shows the values of the Floquet multipliers for different values of $I$. As $I$ tends to $I_6=7.92197768$, the second most negative Floquet multiplier tends to $-1$.

\begin{table}[h!]
\caption{By decreasing the value of I, the multipliers $\mu_4$ decreases, crosses the value -1 for $I= 7.92197768$ and enters into the unit circle.  }\label{table}

\begin{center}
\begin{tabular}{lllll}
\hline\noalign{\smallskip}
I           & $\mu_1$  & $\mu_2$  & $\mu_3$ & $\mu_4$\\
\noalign{\smallskip}\hline\noalign{\smallskip}
    7.92197799  & 1.000 & 0.000  &-2940.687  &-1.041\\
    7.92197793  & 1.000 &-0.000  &-2964.042  &-1.033\\
    7.92197787  & 1.000 & 0.000  &-2987.386  & -1.025\\
    7.92197781  & 1.000 & 0.000  &-3010.719  &-1.017\\
    7.92197775  & 1.000 & -0.000 &-3034.042  &-1.009\\
    7.92197768  & 1.000 & 0.000  & -3057.354 &-1.001\\
    7.92197762  & 1.000 & -0.000 &-3080.655  &-0.993 \\
    7.92197756  & 1.000 & 0.000  &-3103.946  &-0.986\\
    7.92197750  & 1.000 & 0.000  &-3127.225  &-0.978\\
    7.92197743  & 1.000 & 0.000  &-3150.494  &-0.9713\\
\noalign{\smallskip}\hline\noalign{\smallskip}
\end{tabular}
\end{center}
\end{table}

\section{Conclusions}
In 1952 Hodgkin and Huxley developed the pioneer and still
up-to-date mathematical model for describing the activity of the
giant squid axon. Depending on the value of the external current
stimuli, this fourth-order nonlinear dynamical system exhibits many
complex behaviors, such as multiple periodic solutions (both stable
and unstable) and chaos.

Previous works have treated this problem by using several numerical methods, such as shooting and finite difference methods, that are not so simple to handle. In this paper, we jointly exploited shooting, collocation and  harmonic balance methods to obtain the complete bifurcation diagram, therefore detecting all the periodic solutions and the associated bifurcations. In particular, we have shown how the harmonic balance method is extremely handy and works very well in the most complex part of such diagram. Furthermore, harmonic balance and Floquet analysis have permetted to suitably detect the period-doubling bifurcation that entails a route-to-chaos in the HH model. 
% Hassard, Rinzel discovered the different bifurcation in Hodgkin–Huxley model with the variation of overall current I, they used the both numerical methods, shooting for Hassard and other finite difference method for Rinzel, the second method is sufficient for detecting the all stable and unstable periodic solution, but is impractical,(Jacobian matrix in linearized method ...). In our work, we use both methods for detected the unstable solution, the first method is the collocation method, gives good results, but the HB is interesting in our case, a help of the HB, we explored the doubling period bifurcation, it the route to chaos, In future works, the adapt the method , we using the fast-slow dynamic in HH. we use Floquet theory the periodic solutions detected bye HB are highly unstable, consequently, the chaotic solution are also, thus, his study is interesting in the future work.
%In this paper, we have characterized the
%periodic solutions that arise in the HH model through an extremely
%handy collocation method. In future works we envisage to detect the
%periodic solutions of the HH model for all the values of $I$, in
%order to obtain the corresponding bifurcation diagram.

\section*{Acknowledgements}
We would like to thank: R\'{e}gion Haute Normandie, CPER and FEDER (RISC project) for financial
support.

%\section*{References}

\bibliography{HH_bib}

\end{document}